\documentclass[11pt,a4paper]{article}

\usepackage[T1]{fontenc}
\usepackage[utf8]{inputenc}
\usepackage{microtype}

\usepackage{amsmath,amsthm,mathtools}
\usepackage{newtxtext,newtxmath}
\usepackage{graphicx}
\graphicspath{{figures/}}

\usepackage[margin=1in]{geometry}
\usepackage[colorlinks=true,linkcolor=blue,citecolor=blue,urlcolor=blue]{hyperref}

\newtheorem{theorem}{Theorem}[section]
\newtheorem{lemma}[theorem]{Lemma}
\newtheorem{proposition}[theorem]{Proposition}
\newtheorem{corollary}[theorem]{Corollary}
\newtheorem{conjecture}[theorem]{Conjecture}
\theoremstyle{definition}
\newtheorem{definition}[theorem]{Definition}
\theoremstyle{remark}
\newtheorem{remark}[theorem]{Remark}

\newcommand{\R}{\mathbb{R}}
\newcommand{\E}{\mathbb{E}}

\newcommand{\T}{\mathbb{T}}
\newcommand{\Z}{\mathbb{Z}}
\newcommand{\N}{\mathbb{N}}

\newcommand{\p}{\partial}
\newcommand{\fL}{(-\Delta)^{\alpha/2}}
\newcommand{\fLh}{(-\Delta)^{\alpha/4}}
\newcommand{\Capt}{\partial_t^\beta}
\newcommand{\Lie}{\mathcal{L}}
\newcommand{\dd}{\,\mathrm{d}}
\newcommand{\HS}{\mathrm{HS}}
\newcommand{\Tr}{\mathrm{Tr}}
\newcommand{\Leray}{\mathbb{P}}
\newcommand{\la}{\langle}
\newcommand{\ra}{\rangle}
\newcommand{\inner}[2]{\langle #1, #2 \rangle}

\title{Fractional Stochastic Navier--Stokes Equations:\\Local Well-Posedness and Enstrophy Balance}

\author{Joel Saucedo\\
\small Georgia College \& State University\\
\small Department of Chemistry, Physics, \& Astronomy\\
\small Milledgeville, GA 31061, USA}

\date{\today}

\begin{document}
\maketitle

\begin{abstract}
In this paper, we investigate a stochastic incompressible fluid on the three-dimensional torus in which power-law temporal memory, represented by a Caputo derivative of order $\beta\in(1/2,1)$, is coupled to non-local dissipation, represented by a fractional Laplacian of order $\alpha\in(1,2)$, under spatially smooth trace-class noise with a superlinear coefficient. We first derive the equation from constrained Hamiltonian mechanics through a Mori--Zwanzig elimination of fast degrees of freedom, so that the memory kernel and the random force emerge together rather than being postulated; a Tauberian limit identifies the Caputo operator as the exact analytic image of a scale-free harmonic heat bath. Representing vorticity as a differential two-form, we prove that the vortex-stretching functional coincides with the sign-indefinite quadratic form of the strain tensor, and that the velocity is recovered from the vorticity with an exact one-derivative gain. Using a subordination identity for the Mittag--Leffler propagator, we show that the linear kernel is a probability density with an algebraic far-field tail. Within a Gelfand triple adapted to the Stokes operator we establish local existence and uniqueness in $H^s$ for $s>3/2$, with the data and the forcing transported by two distinct propagators, and we prove that the memory threshold $\beta>1/2$ is sharp for the existence of function-valued solutions under white-in-time forcing. We then derive an enstrophy balance, governed by a fractional differential inequality, read off the critical exponent $\sigma=\beta(1+3/\alpha)$ separating the dissipation-dominated regime from possible stretching-driven growth, and close with a Beale--Kato--Majda continuation criterion together with a geometric condition for global existence.
\end{abstract}

\medskip
\noindent\textbf{Keywords:} Navier--Stokes equations; fractional Laplacian; Caputo derivative; stochastic partial differential equations; Gelfand triple; vorticity; enstrophy; Mittag--Leffler function.

\medskip
\noindent\textbf{MSC 2020:} 35Q30, 35R11, 60H15, 76D03, 76M35.

\section{Introduction and Physical Formulation}\label{sec:intro}
% =============================================================================

The fate of an incompressible fluid is decided by an antagonism between two tendencies. Viscosity drains kinetic energy and smooths gradients, while the inertial nonlinearity concentrates vorticity by stretching and folding vortex lines. For the classical three-dimensional Navier--Stokes equations the resolution of this antagonism from smooth data is famously unsettled.\cite{LL6, ConstantinFoias, MajdaBertozzi} We establish here a new structural regime in which the competition between inertial stretching and anomalous dissipation is governed by non-local temporal memory and a fluctuating force. The questions we answer are intrinsic to that regime: when a fluid retains a memory of its past stresses, dissipates through long-range interactions, and is driven by a random force, what is the appropriate evolution equation, how do its geometric and probabilistic structures constrain one another, and what may be proved about its solutions.

Our standpoint is the one familiar from statistical physics, where a macroscopic law is obtained by eliminating fast and microscopic variables from an underlying mechanical description rather than written by analogy.\cite{LL5, Zwanzig} Carried out for a fluid with slow relaxation, this elimination endows the hydrodynamic equation with a memory kernel and a fluctuating force, the two linked by fluctuation and dissipation. Power-law memory yields a Caputo fractional derivative in time, and long-range microscopic interaction yields a fractional Laplacian in space. The nonstandard operators of the model are therefore inherited rather than imposed, and a central purpose of this paper is to make that inheritance precise before any analysis is undertaken.

\subsection{The model}\label{subsec:model}

We study the system
\begin{equation}\label{eq:main}
\begin{aligned}
\Capt \mathbf{u} &= -(\mathbf{u}\cdot\nabla)\mathbf{u} - \nabla p - \nu \fL \mathbf{u} + \Phi(\mathbf{u})\,\dot{W}, \\
\nabla\cdot\mathbf{u} &= 0, \qquad \mathbf{u}|_{t=0} = \mathbf{u}_0,
\end{aligned}
\end{equation}
on the torus $\T^3=(\R/2\pi\Z)^3$, with velocity $\mathbf{u}(t,x)\in\R^3$, pressure $p$, viscosity $\nu>0$, dissipation exponent $\alpha\in(1,2)$, and memory index $\beta\in(0,1)$. Here $\Capt$ is the Caputo derivative of order $\beta$, the operator $\fL$ is the fractional Laplacian, and $\Phi(\mathbf{u})\dot W$ is a stochastic force whose covariance is specified in Section \ref{sec:framework}. The noise is taken to be trace-class, which is the hypothesis under which the second moments remain finite in three dimensions. The periodic setting is chosen because the spectral theory of the Stokes operator is there at its cleanest and the statistical heuristics of Section \ref{sec:scaling} are at their sharpest; the passage to $\R^3$ requires only homogeneous Besov spaces in place of their periodic counterparts.

\subsection{Incompressible flow as constrained Hamiltonian mechanics}\label{subsec:hamiltonian}

Let $H=\{\mathbf{v}\in L^2(\T^3;\R^3):\nabla\cdot\mathbf{v}=0,\ \int_{\T^3}\mathbf{v}\dd x=0\}$ be the configuration space of an incompressible fluid, with kinetic energy
\[
\mathcal{H}[\mathbf{u}]=\tfrac12\int_{\T^3}|\mathbf{u}|^2\dd x.
\]
In the formulation of Arnold,\cite{Arnold1966, ArnoldKhesin} the Euler flow is the geodesic motion generated by $\mathcal{H}$ on the group of volume-preserving diffeomorphisms. Incompressibility is a constraint, and the pressure is the Lagrange multiplier that enforces it. Concretely, the multiplier is realized by the Leray--Helmholtz projector $\Leray$ onto divergence-free fields, introduced by Leray\cite{Leray} and central to the modern theory of Temam,\cite{Temam} which annihilates the gradient part and removes $p$ from the dynamics; in Fourier variables it acts as $\widehat{\Leray f}(\xi)=\hat f(\xi)-|\xi|^{-2}(\xi\cdot\hat f(\xi))\,\xi$, the orthogonal projection of $\hat f(\xi)$ onto the plane perpendicular to $\xi$. This is precisely the mechanism that permits a clean elimination of the pressure in the vorticity formulation of Section \ref{sec:geometry}.

\subsection{Memory and noise from projection}\label{subsec:mz}

A real fluid carries fast degrees of freedom that the hydrodynamic description leaves unresolved. The projection formalism of Mori\cite{Mori} and Zwanzig\cite{Zwanzig} eliminates these exactly, at the price of two new structures in the equation for the slow observables: a memory integral and a fluctuating force. Let $P$ be the orthogonal projection onto the slow modes and let $Q=\mathrm{Id}-P$ be its complement. Applying $P$ to the Liouville equation $\p_t\mathbf{u}=\mathcal L\mathbf{u}$ for the full phase variable and resolving the orthogonal evolution $e^{tQ\mathcal L}$ through Duhamel's identity, one obtains for the resolved velocity the generalized Langevin equation
\begin{equation}\label{eq:gle}
\p_t\mathbf{u}(t)=\Omega\,\mathbf{u}(t)+\int_0^t K(t-s)\,\mathbf{u}(s)\dd s+\eta(t),
\end{equation}
in which the streaming term is $\Omega=P\mathcal L$, the memory kernel is $K(t)=P\mathcal L\,e^{tQ\mathcal L}\,Q\mathcal L$, and the residual force is $\eta(t)=e^{tQ\mathcal L}Q\mathcal L\,\mathbf{u}(0)$. Kernel and force are not independent: linear-response theory ties them through the stationary autocorrelation of the force by the second fluctuation--dissipation relation,
\begin{equation}\label{eq:fdt}
K(t)=\frac{\la\eta(t),\eta(0)\ra}{\la\mathbf{u},\mathbf{u}\ra},
\end{equation}
which is the content of the Kubo theory of dissipation.\cite{Kubo} Two regimes of \eqref{eq:gle} concern us. When $K$ is short-ranged, scaling its support to zero collapses the memory integral into the instantaneous friction $\nu\,\p_t\mathbf{u}$ and reproduces ordinary viscous dissipation.

The power-law regime, by contrast, is selected by a single structural hypothesis on the bath of unresolved modes. Following the standard oscillator realization of the generalized Langevin equation,\cite{Zwanzig} model the bath by a continuum of harmonic oscillators coupled bilinearly to the slow variable, with spectral density $J(\omega)$. The memory kernel is then the cosine transform of the friction spectrum,
\begin{equation}\label{eq:cosinetransform}
K(t)=\frac{2}{\pi}\int_0^\infty \frac{J(\omega)}{\omega}\cos(\omega t)\dd\omega .
\end{equation}
Call the bath \emph{scale-free} when its low-frequency spectrum is self-similar,
\begin{equation}\label{eq:scalefree}
J(\omega)\sim c\,\omega^{\beta}\qquad(\omega\to0),\qquad \beta\in(0,1).
\end{equation}
The large-$t$ behavior of \eqref{eq:cosinetransform} is governed by the $\omega\to0$ behavior of its integrand. Substituting \eqref{eq:scalefree} gives $J(\omega)/\omega\sim c\,\omega^{\beta-1}$, and the elementary cosine transform
\begin{equation}\label{eq:cosineeval}
\int_0^\infty \omega^{\beta-1}\cos(\omega t)\dd\omega=\Gamma(\beta)\,\cos\frac{\pi\beta}{2}\;t^{-\beta}\qquad(0<\beta<1),
\end{equation}
combined with the Abelian theorem for regularly varying transforms of Bingham, Goldie, and Teugels,\cite{Bingham} yields the algebraically decaying kernel
\begin{equation}\label{eq:powerlawkernel}
K(t)\sim \frac{C_\beta}{t^{\beta}}\qquad(t\to\infty),\qquad C_\beta=\frac{2c}{\pi}\,\Gamma(\beta)\,\cos\frac{\pi\beta}{2}.
\end{equation}
This power law is the only input we take from the microscopic picture. Inserting \eqref{eq:powerlawkernel} into the memory integral of \eqref{eq:gle}, the convolution against $t^{-\beta}$ is by definition a Riemann--Liouville fractional integral,
\begin{equation}\label{eq:RLid}
\int_0^t \frac{\mathbf{u}(s)}{(t-s)^{\beta}}\dd s=\Gamma(1-\beta)\,I^{1-\beta}\mathbf{u}(t),
\end{equation}
and differentiation in time turns the resulting term into the Caputo derivative, because $\p_t I^{1-\beta}$ acting on the velocity increment is precisely $\Capt$; this is the defining relation of the Caputo operator (Podlubny;\cite{Podlubny} Kilbas, Srivastava, and Trujillo\cite{Kilbas}). Thus the hypothesis of a self-similar bath, and nothing more, converts the memory term into $\Capt$. When in addition the eliminated interactions are long-ranged in space, with a likewise scale-free spatial spectrum, the effective spatial generator is the fractional Laplacian $\fL$, whose Fourier symbol $|\xi|^\alpha$ generates the isotropic $\alpha$-stable L\'evy semigroup and so encodes anomalous, jump-type transport (Caffarelli and Silvestre\cite{CaffarelliSilvestre}). The Gaussian force $\eta$, arising through the central-limit structure of the sum over eliminated modes, is the term $\Phi(\mathbf{u})\dot W$, and its covariance is fixed by the same kernel \eqref{eq:powerlawkernel} through the relation \eqref{eq:fdt}.

The three nonstandard ingredients of \eqref{eq:main} are therefore not independent modeling choices but the generic residue of eliminating fast degrees of freedom from a Hamiltonian fluid with a scale-free bath, so that probability enters as a derived rather than a primitive feature of the description. We adopt the self-similarity \eqref{eq:scalefree} as the one axiomatic assumption on the unresolved dynamics; everything downstream is its consequence.

\subsection{Emergence of Caputo calculus}\label{subsec:caputo-emergence}

The passage from the power-law memory of \eqref{eq:gle} to the Caputo operator of \eqref{eq:main} was presented above as a modeling step. We now record it as precise mathematics, so that the fractional derivative is the rigorous analytic image of a microscopic memory rather than a postulate. The argument has three layers: an explicit mechanical model whose elimination produces the generalized Langevin equation \eqref{eq:gle}, a Tauberian theorem converting the scale-free hypothesis \eqref{eq:scalefree} into the power-law kernel \eqref{eq:powerlawkernel}, and an exact operator identity converting that kernel into $\Capt$.

\begin{lemma}[Harmonic heat-bath realization]\label{lem:bath}
Let a slow scalar mode $q(t)$ be coupled to a bath of $N$ harmonic oscillators through the Ford--Kac--Zwanzig Hamiltonian
\[
H_N=\tfrac12\dot q^2+\sum_{k=1}^N\Big[\tfrac12\dot x_k^2+\tfrac12\omega_k^2\Big(x_k-\frac{g_k}{\omega_k^2}\,q\Big)^2\Big].
\]
Then elimination of the bath variables yields the exact generalized Langevin equation
\[
\ddot q(t)=-\int_0^t K_N(t-s)\,\dot q(s)\dd s+\eta_N(t),\qquad
K_N(t)=\sum_{k=1}^N\frac{g_k^2}{\omega_k^2}\cos(\omega_k t),
\]
with $\eta_N$ a stationary Gaussian process whose covariance is $\la\eta_N(t),\eta_N(0)\ra=k_BT\,K_N(t)$ in thermal equilibrium. If the couplings are distributed so that $\tfrac{\pi}{2}\sum_k(g_k^2/\omega_k^2)\,\delta(\omega-\omega_k)\rightharpoonup J(\omega)/\omega$ as $N\to\infty$, then $K_N\to K$ of \eqref{eq:cosinetransform} and the force covariance converges to $K$, which is the relation \eqref{eq:fdt}.
\end{lemma}

\begin{proof}
Each bath coordinate obeys the forced linear equation $\ddot x_k+\omega_k^2 x_k=g_k q$, whose solution by variation of constants is
\[
x_k(t)=x_k(0)\cos\omega_k t+\frac{\dot x_k(0)}{\omega_k}\sin\omega_k t+\frac{g_k}{\omega_k}\int_0^t\sin\!\big(\omega_k(t-s)\big)q(s)\dd s.
\]
Inserting this into the equation $\ddot q=-\sum_k g_k(x_k-g_k q/\omega_k^2)$ for the slow mode and integrating the memory term once by parts converts the $\sin$ kernel into the $\cos$ kernel $K_N(t)=\sum_k(g_k^2/\omega_k^2)\cos(\omega_k t)$ acting on $\dot q$, while the terms carrying the initial bath data assemble into $\eta_N$. Drawing $x_k(0),\dot x_k(0)$ from the Gibbs measure of $H_N$ gives $\la\eta_N(t),\eta_N(0)\ra=k_BT\,K_N(t)$, the second fluctuation--dissipation relation. The stated weak convergence of the coupling measure turns the sum defining $K_N$ into the integral \eqref{eq:cosinetransform}; this is the Ford--Kac construction.\cite{Zwanzig, Kubo}
\end{proof}

\begin{theorem}[Tauberian emergence of the fractional memory]\label{thm:tauberian}
Suppose the bath spectral density satisfies the scale-free law \eqref{eq:scalefree}, $J(\omega)\sim c\,\omega^{\beta}$ as $\omega\to0^+$ with $\beta\in(0,1)$, and that $J(\omega)/\omega$ is locally integrable and ultimately monotone. Then the kernel \eqref{eq:cosinetransform} is regularly varying of index $-\beta$ at infinity, with
\[
K(t)=C_\beta\,t^{-\beta}\big(1+o(1)\big)\qquad(t\to\infty),\qquad C_\beta=\frac{2c}{\pi}\,\Gamma(\beta)\cos\frac{\pi\beta}{2}>0.
\]
\end{theorem}

\begin{proof}
Fix $\delta>0$ and split $K(t)=\tfrac{2}{\pi}\big(\int_0^\delta+\int_\delta^\infty\big)\,(J(\omega)/\omega)\cos(\omega t)\dd\omega$. On the high-frequency part the integrand is integrable and oscillatory, so by the Riemann--Lebesgue lemma $\int_\delta^\infty(J(\omega)/\omega)\cos(\omega t)\dd\omega=o(t^{-\beta})$ as $t\to\infty$. On the low-frequency part, \eqref{eq:scalefree} gives $J(\omega)/\omega=c\,\omega^{\beta-1}(1+o(1))$, and the exact cosine transform \eqref{eq:cosineeval},
\[
\frac{2}{\pi}\int_0^\infty c\,\omega^{\beta-1}\cos(\omega t)\dd\omega=\frac{2c}{\pi}\,\Gamma(\beta)\cos\frac{\pi\beta}{2}\;t^{-\beta},
\]
supplies the leading term; the monotone-density form of the Abelian--Tauberian theorem of Bingham, Goldie, and Teugels\cite{Bingham} controls the remainder as $o(t^{-\beta})$. Collecting the two contributions yields the stated asymptotic, which is exactly \eqref{eq:powerlawkernel}.
\end{proof}

\begin{theorem}[Kernel--Caputo identity]\label{thm:caputo-id}
Let $X$ be a Banach space and $u\in C^1([0,T];X)$. Then for every $\beta\in(0,1)$,
\begin{equation}\label{eq:caputo-id}
\frac{1}{\Gamma(1-\beta)}\,\frac{\dd}{\dd t}\int_0^t (t-s)^{-\beta}\,u(s)\dd s
=\Capt u(t)+\frac{t^{-\beta}}{\Gamma(1-\beta)}\,u(0).
\end{equation}
In particular, the once-differentiated power-law memory functional with kernel \eqref{eq:powerlawkernel}, acting on a history that starts from rest, $u(0)=0$, coincides exactly with the Caputo derivative $\Capt u$.
\end{theorem}

\begin{proof}
Substitute $r=t-s$ in the memory integral, so that
\[
\int_0^t (t-s)^{-\beta}u(s)\dd s=\int_0^t r^{-\beta}\,u(t-r)\dd r .
\]
The variable $t$ now appears both in the upper limit and in the integrand. Differentiating, and using that the kernel singularity $r^{-\beta}$ is integrable for $\beta<1$ so that the Leibniz rule applies,
\[
\frac{\dd}{\dd t}\int_0^t r^{-\beta}u(t-r)\dd r
= t^{-\beta}\,u(0)+\int_0^t r^{-\beta}\,u'(t-r)\dd r .
\]
Undoing the substitution in the remaining integral gives $\int_0^t (t-s)^{-\beta}u'(s)\dd s=\Gamma(1-\beta)\,\Capt u(t)$, by the definition of the Caputo derivative. Dividing by $\Gamma(1-\beta)$ produces \eqref{eq:caputo-id}. When $u(0)=0$ the transient $t^{-\beta}u(0)/\Gamma(1-\beta)$ vanishes and the left side reduces to $\Capt u(t)$ (Podlubny;\cite{Podlubny} Kilbas, Srivastava, and Trujillo\cite{Kilbas}).
\end{proof}

Lemma \ref{lem:bath} together with Theorems \ref{thm:tauberian} and \ref{thm:caputo-id} closes the chain: a scale-free harmonic bath produces, in the continuum limit, the power-law memory kernel, and the once-differentiated memory operator is exactly $\Capt$. The fractional Laplacian arises by the same Tauberian principle applied to a scale-free \emph{spatial} spectrum, where the regularly varying symbol $|\xi|^\alpha$ is the generator of the isotropic $\alpha$-stable L\'evy semigroup (Caffarelli and Silvestre;\cite{CaffarelliSilvestre} Meerschaert and Sikorskii\cite{MeerschaertSikorskii}). The operators of \eqref{eq:main} are thus derived, not imposed.

\subsection{Main results and outline}\label{subsec:results}

The emergence of the Caputo operator is itself made rigorous in Section \ref{subsec:caputo-emergence}. Section \ref{sec:framework} assembles the analytic framework: Littlewood--Paley theory, the fractional operators, the two Mittag-Leffler propagators of the linear evolution, and the Gelfand triple in which the trace-class noise lives. Section \ref{sec:geometry} develops the geometry of vorticity, proves the sign-indefinite stretching identity and the exact Biot--Savart gain, and establishes the algebraic far-field asymptotics of the linear kernel by subordination. Section \ref{sec:analysis} contrasts the fractional and classical evolutions, proves local well-posedness of the vorticity equation in $H^s$ for $s>3/2$, derives the enstrophy balance through a fractional Leibniz inequality, and gives a continuation criterion with a geometric condition for global existence. Section \ref{sec:scaling} extracts the scaling laws of the critical exponents and indicates the statistical and conditional directions that the framework opens.

The principal results are: the identification of the Caputo operator as the exact analytic limit of a scale-free harmonic heat bath (Lemma \ref{lem:bath} and Theorems \ref{thm:tauberian} and \ref{thm:caputo-id}); the explicit separation of the data and forcing propagators (Proposition \ref{prop:propagator}); the sign-indefinite character of vortex stretching as a consequence of its differential-form structure (Proposition \ref{prop:stretch}); the probability-density and heavy-tail structure of the linear kernel (Proposition \ref{prop:kernel}); the sharp memory threshold $\beta>\tfrac12$ for the existence of function-valued solutions (Proposition \ref{prop:beta-sharp}); local existence and uniqueness with trace-class noise (Theorem \ref{thm:lwp}); the enstrophy balance (Theorem \ref{thm:enstrophy}) with its scaling corollaries; and a Beale--Kato--Majda continuation criterion with conditional global existence (Theorem \ref{thm:continuation} and Corollary \ref{cor:global}).

\begin{remark}[Global behavior and superlinear noise]\label{rem:mijena}
The superlinear coefficient $\Phi(\mathbf{u})\sim|\mathbf{u}|^{1+\gamma}$ is precisely the structure under which time-fractional stochastic equations are known to lose global solutions. For the time-fractional stochastic heat equation, Mijena and Nane\cite{MijenaNane} constructed the solution theory and analyzed its moment growth, and Asogwa, Mijena, and Nane\cite{AsogwaMijenaNane} proved finite-time blow-up of moments under superlinear forcing; related asymptotics appear in Foondun and Nane.\cite{FoondunNane} We therefore do not expect global existence in the superlinear regime, which is the reason our well-posedness statement is local, and the enstrophy balance of Section \ref{sec:analysis} is the natural diagnostic for the onset of growth.
\end{remark}

% =============================================================================
\section{Analytical Framework}\label{sec:framework}
% =============================================================================

\subsection{Littlewood--Paley decomposition and Besov spaces}\label{subsec:lp}

Following Bahouri, Chemin, and Danchin,\cite{BCD} let $\{\Delta_j\}_{j\ge-1}$ be the standard nonhomogeneous Littlewood--Paley projections associated with a dyadic partition of unity, so that $\Delta_j$ localizes frequencies to the shell $|\xi|\sim2^j$. For $s\in\R$ and $p,q\in[1,\infty]$ the Besov space $B^s_{p,q}$ consists of the tempered distributions with finite norm
\[
\|f\|_{B^s_{p,q}}:=\Big(\sum_{j\ge-1}2^{jsq}\|\Delta_j f\|_{L^p}^q\Big)^{1/q},
\]
modified in the usual way for $q=\infty$ (see also Triebel\cite{Triebel}). We abbreviate $H^s=B^s_{2,2}$, which carries the equivalent Fourier norm $\|f\|_{H^s}^2=\sum_{\xi\in\Z^3}\la\xi\ra^{2s}|\hat f(\xi)|^2$ with $\la\xi\ra=(1+|\xi|^2)^{1/2}$. Two facts, valid for $s>3/2$ on $\T^3$ and established in the same references, are used repeatedly. The first is the embedding and Banach-algebra pair
\begin{equation}\label{eq:algebra}
H^s\hookrightarrow L^\infty,\qquad \|fg\|_{H^s}\le C_s\|f\|_{H^s}\|g\|_{H^s};
\end{equation}
the embedding holds because $\sum_\xi|\hat f(\xi)|\le(\sum_\xi\la\xi\ra^{-2s})^{1/2}\|f\|_{H^s}$ is finite exactly when $2s>3$, and the product law follows from it by frequency decomposition. The second is the tame multiplication law
\begin{equation}\label{eq:tame}
\|fg\|_{H^{s-1}}\le C_s\|f\|_{H^s}\|g\|_{H^{s-1}}\qquad(s>3/2),
\end{equation}
which keeps the lower-order factor at its own regularity. The threshold $s>3/2$ recurs throughout precisely because the superlinear coefficient $\Phi$ requires control of the velocity in $L^\infty$, which \eqref{eq:algebra} supplies.

\subsection{Fractional operators and the linear propagators}\label{subsec:frac}

The fractional Laplacian $\fL$ is the Fourier multiplier with symbol $|\xi|^\alpha$, acting on the Fourier coefficients indexed by $\xi\in\Z^3$. Writing $I^\gamma g(t)=\Gamma(\gamma)^{-1}\int_0^t(t-s)^{\gamma-1}g(s)\dd s$ for the Riemann--Liouville integral of order $\gamma>0$, the Caputo derivative of order $\beta\in(0,1)$ of a function $f\in C^1([0,T])$ is the once-integrated derivative,
\[
\Capt f(t)=\frac{1}{\Gamma(1-\beta)}\int_0^t(t-s)^{-\beta}f'(s)\dd s=I^{1-\beta}f'(t),
\]
in the sense of Caputo as developed by Podlubny,\cite{Podlubny} by Kilbas, Srivastava, and Trujillo,\cite{Kilbas} and by Diethelm.\cite{Diethelm} The linear evolution is governed by the two-parameter Mittag-Leffler function $E_{a,b}(z)=\sum_{k\ge0}z^k/\Gamma(ak+b)$, the natural eigenfunction of $\Capt$. From the monograph of Gorenflo, Kilbas, Mainardi, and Rogosin\cite{Gorenflo} we take, for $\beta\in(0,1)$ and $x\ge0$, the two-sided bounds
\begin{equation}\label{eq:MLbounds}
0\le E_{\beta,1}(-x)\le 1,\qquad 0\le E_{\beta,\beta}(-x)\le \frac{C_\beta}{1+x},
\end{equation}
in which both functions are completely monotone in $x$, together with the large-argument asymptotic
\begin{equation}\label{eq:MLasym}
E_{\beta,1}(-x)=\frac{1}{\Gamma(1-\beta)}\,\frac{1}{x}+O(x^{-2})\qquad(x\to\infty),
\end{equation}
which exhibits the algebraic, rather than exponential, relaxation characteristic of fractional dynamics.

\begin{proposition}[Linear propagators]\label{prop:propagator}
The solution of $\Capt u=-\nu\fL u$ with $u(0)=u_0$ is $u(t)=S_\beta(t)u_0$, where $S_\beta(t)$ is the Fourier multiplier with symbol $E_{\beta,1}(-\nu t^\beta|\xi|^\alpha)$. The solution of $\Capt u=-\nu\fL u+f$ with $u(0)=0$ is
\[
u(t)=\int_0^t(t-s)^{\beta-1}P_\beta(t-s)\,f(s)\dd s,\qquad \widehat{P_\beta(\tau)}(\xi)=E_{\beta,\beta}(-\nu\tau^\beta|\xi|^\alpha).
\]
The data are transported by $E_{\beta,1}$ and the forcing by the distinct kernel $(t-s)^{\beta-1}E_{\beta,\beta}$.
\end{proposition}

\begin{proof}
Expand $u$ in a Fourier series; for each $\xi\in\Z^3$ the equation \eqref{eq:main} reduces, after projection, to the scalar Caputo problem
\[
\Capt\hat u(t,\xi)=-\nu|\xi|^\alpha\hat u(t,\xi)+\hat f(t,\xi),\qquad \hat u(0,\xi)=\hat u_0(\xi).
\]
Apply the Laplace transform in time, writing $\tilde u(z,\xi)=\int_0^\infty e^{-zt}\hat u(t,\xi)\dd t$. The Caputo derivative transforms, by Podlubny,\cite{Podlubny} as $\mathcal{L}\{\Capt g\}(z)=z^\beta\tilde g(z)-z^{\beta-1}g(0)$, so the transformed equation reads
\[
z^\beta\tilde u-z^{\beta-1}\hat u_0+\nu|\xi|^\alpha\tilde u=\tilde f,
\]
whence
\[
\tilde u(z,\xi)=\frac{z^{\beta-1}}{z^\beta+\nu|\xi|^\alpha}\,\hat u_0(\xi)+\frac{1}{z^\beta+\nu|\xi|^\alpha}\,\tilde f(z,\xi).
\]
The two elementary Mittag-Leffler Laplace identities recorded by Gorenflo, Kilbas, Mainardi, and Rogosin,\cite{Gorenflo}
\[
\mathcal{L}^{-1}\!\Big\{\frac{z^{\beta-1}}{z^\beta+a}\Big\}(t)=E_{\beta,1}(-at^\beta),
\qquad
\mathcal{L}^{-1}\!\Big\{\frac{1}{z^\beta+a}\Big\}(t)=t^{\beta-1}E_{\beta,\beta}(-at^\beta),
\]
specialized to $a=\nu|\xi|^\alpha$, invert the first term to $E_{\beta,1}(-\nu t^\beta|\xi|^\alpha)\hat u_0(\xi)$ and, by the convolution theorem applied to the second, give
\[
\hat u(t,\xi)=E_{\beta,1}(-\nu t^\beta|\xi|^\alpha)\hat u_0(\xi)+\int_0^t (t-s)^{\beta-1}E_{\beta,\beta}(-\nu(t-s)^\beta|\xi|^\alpha)\hat f(s,\xi)\dd s.
\]
Summing over $\xi$ yields the stated formula. The interchange of summation and Laplace inversion is justified by \eqref{eq:MLbounds}, which makes each symbol bounded and the series absolutely convergent in $H^s$.
\end{proof}

\begin{remark}\label{rem:twoprop}
The factor $(t-s)^{\beta-1}$ multiplying the forcing propagator is integrable at $s=t$ for $\beta\in(0,1)$ and is responsible for the singularity structure of the stochastic convolution; the data propagator carries no such factor. The two kernels must be kept distinct, since identifying them would alter every time exponent in Section \ref{sec:analysis}.
\end{remark}

\begin{figure}[t]
\centering
\includegraphics[width=0.9\textwidth]{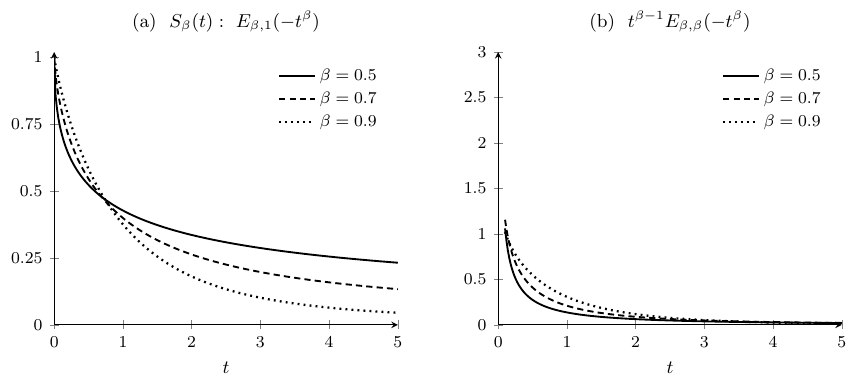}
\caption{Representative profiles of the two propagators of Proposition \ref{prop:propagator}, for $\nu=1$ and $|\xi|=1$. (a) The datum is transported by the Mittag-Leffler function $E_{\beta,1}(-t^\beta)$, which decays monotonically with an algebraic tail that is heavier for smaller $\beta$. (b) The forcing enters through the distinct, singular Duhamel kernel $t^{\beta-1}E_{\beta,\beta}(-t^\beta)$, which is integrable at $t=0$ for $\beta\in(0,1)$. Identifying the two kernels would alter the time exponents of Section \ref{sec:analysis}; see Remark \ref{rem:twoprop}.}
\label{fig:propagators}
\end{figure}

\subsection{Gelfand triple and trace-class noise}\label{subsec:gelfand}

White noise is a generalized function and is not realized in $L^2$; its proper setting is the rigged Hilbert space of Gelfand and Vilenkin.\cite{Gelfand} Let $A=\Leray\fL$ be the fractional Stokes operator on $H$ with domain $D(A)$, set $V=D(A^{1/2})$, and let $V'$ be its dual. Then
\[
V\hookrightarrow H\hookrightarrow V'
\]
is a Gelfand triple with dense and continuous embeddings, the standard analytic frame for variational solutions of stochastic evolution equations (Prévôt and Röckner;\cite{PrevotRockner} Liu and Röckner\cite{LiuRockner}). Let $\{e_k\}_{k\in\N}$ be the orthonormal eigenbasis of $A$ on $H$, which exists because $A$ has compact resolvent on the torus.

\begin{definition}[$Q$-Wiener process; spatially smooth noise]\label{def:qwiener}
Let $Q:H\to H$ be symmetric, non-negative, with eigenvalues $q_k\ge0$ in the basis $\{e_k\}$. Following Da~Prato and Zabczyk,\cite{DaPratoZabczyk} the associated $Q$-Wiener process is $W(t)=\sum_k\sqrt{q_k}\,\beta_k(t)\,e_k$ with independent standard Brownian motions $\beta_k$, and its covariance is $\E[\inner{W(t)}{a}\inner{W(s)}{b}]=(t\wedge s)\inner{Qa}{b}$. We call the noise \emph{spatially smooth} of order $r$ if
\[
\|Q^{1/2}\|_{\HS(H,H^{r})}^2=\sum_{k}q_k\,\|e_k\|_{H^{r}}^2<\infty.
\]
We assume throughout that the noise is spatially smooth of order $s+1$, that is
\[
\|Q^{1/2}\|_{\HS(H,H^{s+1})}<\infty.
\]
\end{definition}

For predictable integrands the It\^o isometry reads
\begin{equation}\label{eq:ito}
\E\Big\|\int_0^t\Psi(\sigma)\dd W(\sigma)\Big\|_{H^s}^2=\int_0^t\E\|\Psi(\sigma)Q^{1/2}\|_{\HS(H,H^s)}^2\dd\sigma.
\end{equation}
The requirement $\Tr Q<\infty$ is essential: for space-time white noise $Q=\mathrm{Id}$ has infinite trace in three dimensions, the stochastic convolution then has infinite $H$-norm, and pointwise second moments diverge. Coloring the noise so that $Q$ is trace-class, and smooth of order $s+1$, restores finite second moments and is the standard hypothesis for function-valued solutions in three dimensions (Da~Prato and Zabczyk;\cite{DaPratoZabczyk} Flandoli\cite{Flandoli}).

The diffusion coefficient is obtained by multiplying with a superlinear scalar and projecting the result back onto divergence-free fields, so that the stochastic force remains in $H$. Precisely, fix $\gamma\in(0,1/2)$ and define $\Phi(\mathbf{u}):H\to H$ by
\begin{equation}\label{eq:phidef}
\Phi(\mathbf{u})\,h=\Leray\big(|\mathbf{u}|^{1+\gamma}\,h\big),\qquad h\in H,
\end{equation}
so that the stochastic term in \eqref{eq:main} is the divergence-free field $\Phi(\mathbf{u})\dd W=\Leray\big(|\mathbf{u}|^{1+\gamma}\dd W\big)$, with the projector applied last. We record the mapping property used throughout Section \ref{sec:analysis}. For $s>3/2$ the space $H^s$ is a Banach algebra, the Nemytskii map $\mathbf{u}\mapsto|\mathbf{u}|^{1+\gamma}$ sends $H^s$ into $H^s$ with $\||\mathbf{u}|^{1+\gamma}\|_{H^s}\le C(1+\|\mathbf{u}\|_{H^s})^{1+\gamma}$, and $\Leray$ is bounded on every $H^r$. Hence for any orthonormal basis $\{e_k\}$,
\begin{equation}\label{eq:phiHS}
\begin{aligned}
\|\Phi(\mathbf{u})Q^{1/2}\|_{\HS(H,H^s)}^2
&=\sum_k q_k\big\|\Leray(|\mathbf{u}|^{1+\gamma}e_k)\big\|_{H^s}^2
\le C\,\||\mathbf{u}|^{1+\gamma}\|_{H^s}^2\sum_k q_k\|e_k\|_{H^s}^2\\
&\le C\,(1+\|\mathbf{u}\|_{H^s})^{2(1+\gamma)}\,\|Q^{1/2}\|_{\HS(H,H^s)}^2,
\end{aligned}
\end{equation}
which is finite, since smoothness of order $s+1$ implies smoothness of order $s$. Thus $\Phi(\mathbf{u})Q^{1/2}$ is a Hilbert--Schmidt map $H\to H^s$, the It\^o isometry \eqref{eq:ito} applies to it, and \eqref{eq:phidef} realizes the superlinear regime of Remark \ref{rem:mijena}.

% =============================================================================
\section{Vorticity Geometry and Linear Kernel Estimates}\label{sec:geometry}
% =============================================================================

\subsection{Vorticity as a two-form and the sign of stretching}\label{subsec:twoform}

Let $\omega=\nabla\times\mathbf{u}$. In the language of differential forms, used here in the manner of Arnold and Khesin,\cite{ArnoldKhesin} the velocity one-form $u^\flat=g(\mathbf{u},\cdot)$ has exterior derivative $\dd u^\flat$, the two-form that is the geometric incarnation of vorticity, and transport along the flow is governed by the Lie derivative through Cartan's magic formula $\Lie_{\mathbf u}=\dd\,\iota_{\mathbf u}+\iota_{\mathbf u}\,\dd$. Taking the curl of \eqref{eq:main} removes the pressure, since $\nabla\times\nabla p=0$, and gives the vorticity equation in the form classical for incompressible flow (Majda and Bertozzi\cite{MajdaBertozzi}),
\begin{equation}\label{eq:vorticity}
\Capt\omega=-\nu\fL\omega-(\mathbf{u}\cdot\nabla)\omega+(\omega\cdot\nabla)\mathbf{u}+\nabla\times\big(\Phi(\mathbf{u})\dot W\big).
\end{equation}
The advection term transports vorticity and is conservative in $L^2$; the term $(\omega\cdot\nabla)\mathbf{u}$ is vortex stretching, the sole source of nonlinear amplification. Its character is fixed by the following identity.

\begin{proposition}[Stretching identity]\label{prop:stretch}
Let $\mathbf{u}$ be a smooth divergence-free field on $\T^3$ with vorticity $\omega$, and let $S=\tfrac12(\nabla\mathbf{u}+(\nabla\mathbf{u})^{\mathsf T})$ be the strain-rate tensor, with eigenvalues $\lambda_1\ge\lambda_2\ge\lambda_3$ and orthonormal eigenframe $\{f_1,f_2,f_3\}$. Then
\[
\inner{(\omega\cdot\nabla)\mathbf{u}}{\omega}_{L^2}=\int_{\T^3}\omega^{\mathsf T}S\,\omega\dd x=\int_{\T^3}\sum_{i=1}^3\lambda_i\,(\omega\cdot f_i)^2\dd x.
\]
Since $\Tr S=\nabla\cdot\mathbf{u}=0$ one has $\lambda_1+\lambda_2+\lambda_3=0$, hence $\lambda_1\ge0\ge\lambda_3$, and the integrand is sign-indefinite: it is positive where $\omega$ aligns with $f_1$ and negative where $\omega$ aligns with $f_3$, unless $S\equiv0$.
\end{proposition}

\begin{proof}
Writing components, $((\omega\cdot\nabla)\mathbf{u})_i=\omega_k\p_k u_i$, so
\[
\inner{(\omega\cdot\nabla)\mathbf{u}}{\omega}_{L^2}=\int_{\T^3}\omega_k(\p_k u_i)\omega_i\dd x=\int_{\T^3}\omega_i\omega_k\,\p_k u_i\dd x.
\]
Decompose the velocity gradient into symmetric and antisymmetric parts, $\p_k u_i=S_{ik}+A_{ik}$ with $S_{ik}=\tfrac12(\p_k u_i+\p_i u_k)$ and $A_{ik}=\tfrac12(\p_k u_i-\p_i u_k)$. The product $\omega_i\omega_k$ is symmetric under $i\leftrightarrow k$, while $A_{ik}$ is antisymmetric, so their contraction vanishes:
\[
\sum_{i,k}\omega_i\omega_k A_{ik}=0.
\]
There remains $\int\sum_{i,k}\omega_i\omega_k S_{ik}\dd x=\int\omega^{\mathsf T}S\omega\dd x$. Diagonalizing $S=\sum_i\lambda_i f_if_i^{\mathsf T}$ pointwise gives $\omega^{\mathsf T}S\omega=\sum_i\lambda_i(\omega\cdot f_i)^2$. Finally $\Tr S=\p_iu_i=\nabla\cdot\mathbf u=0$ forces $\sum_i\lambda_i=0$; if $S\not\equiv0$ then $\lambda_1>0>\lambda_3$ on a set of positive measure, and the quadratic form takes both signs as the orientation of $\omega$ varies.
\end{proof}

\begin{remark}[No coercive lower bound]\label{rem:nocoercive}
A bound $\inner{(\omega\cdot\nabla)\mathbf{u}}{\omega}_{L^2}\ge C\|\omega\|_X^2$ is excluded on two independent grounds. First, the left side is cubic in $\omega$: since $\mathbf{u}$ depends linearly on $\omega$ through the Biot--Savart law of Proposition \ref{prop:bsgain}, the substitution $\omega\mapsto\lambda\omega$ sends the left side to $\lambda^3$ times its value and the right side to $\lambda^2$ times its value, so no inequality can hold for all $\lambda>0$. Second, by Proposition \ref{prop:stretch} the left side is not of one sign. Consequently the natural scalar object of Section \ref{sec:analysis} is an exact balance and not a one-sided inequality.
\end{remark}

\begin{figure}[t]
\centering
\includegraphics[width=0.74\textwidth]{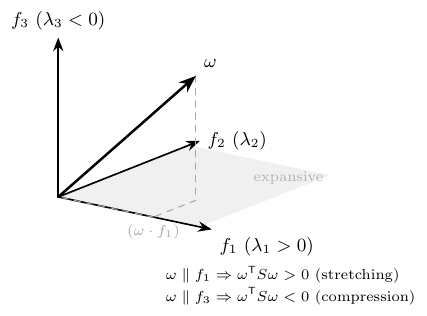}
\caption{The vortex-stretching density in the strain eigenframe of Proposition \ref{prop:stretch}. The strain tensor $S$ has orthonormal eigenframe $\{f_1,f_2,f_3\}$ with eigenvalues summing to zero, so the trace-free constraint forces $\lambda_1\ge0\ge\lambda_3$. The pointwise density is the quadratic form $\omega^{\mathsf T}S\,\omega=\sum_i\lambda_i(\omega\cdot f_i)^2$, which is positive when $\omega$ aligns with the expansive direction $f_1$ and negative when it aligns with the compressive direction $f_3$. The form therefore admits no universal coercive lower bound, which is why the analysis of Section \ref{sec:analysis} rests on an exact enstrophy balance rather than a one-sided inequality.}
\label{fig:strain}
\end{figure}

\subsection{Fractional Biot--Savart and kernel asymptotics}\label{subsec:bs}

\begin{proposition}[Exact Biot--Savart gain]\label{prop:bsgain}
For every mean-zero divergence-free field $\mathbf{u}$ on $\T^3$ with vorticity $\omega=\nabla\times\mathbf{u}$ and every $s\in\R$,
\[
\|\mathbf{u}\|_{H^{s+1}}\le C\,\|\omega\|_{H^s}.
\]
\end{proposition}

\begin{proof}
In Fourier variables divergence-freeness reads $\xi\cdot\hat{\mathbf u}(\xi)=0$, and $\hat\omega(\xi)=i\xi\times\hat{\mathbf u}(\xi)$. Taking the cross product with $\xi$ and using the vector identity $\xi\times(\xi\times\hat{\mathbf u})=(\xi\cdot\hat{\mathbf u})\xi-|\xi|^2\hat{\mathbf u}=-|\xi|^2\hat{\mathbf u}$ gives
\[
i\xi\times\hat\omega(\xi)=i\xi\times(i\xi\times\hat{\mathbf u})=-\,\xi\times(\xi\times\hat{\mathbf u})=|\xi|^2\hat{\mathbf u}(\xi),
\]
so $\hat{\mathbf u}(\xi)=i\xi\times\hat\omega(\xi)/|\xi|^2$ for $\xi\ne0$. Because $\hat{\mathbf u}\perp\xi$ we also have $|\hat\omega(\xi)|=|\xi\times\hat{\mathbf u}(\xi)|=|\xi|\,|\hat{\mathbf u}(\xi)|$, whence the exact relation
\[
|\hat{\mathbf u}(\xi)|=|\xi|^{-1}|\hat\omega(\xi)|,\qquad \xi\ne0.
\]
Therefore
\[
\|\mathbf u\|_{H^{s+1}}^2=\sum_{\xi\ne0}\la\xi\ra^{2(s+1)}|\hat{\mathbf u}(\xi)|^2=\sum_{\xi\ne0}\la\xi\ra^{2(s+1)}|\xi|^{-2}|\hat\omega(\xi)|^2\le C\sum_{\xi\ne0}\la\xi\ra^{2s}|\hat\omega(\xi)|^2=C\|\omega\|_{H^s}^2,
\]
using $\la\xi\ra^{2}|\xi|^{-2}\le C$ for $|\xi|\ge1$ on $\Z^3\setminus\{0\}$.
\end{proof}

The decay of the linear kernel is most transparent through subordination, which expresses the time-fractional propagator as an average of the $\alpha$-stable semigroup over a random time change, in the manner systematized by Meerschaert and Sikorskii.\cite{MeerschaertSikorskii}

\begin{proposition}[Linear kernel as a heavy-tailed probability density]\label{prop:kernel}
Let $G_t$ be the convolution kernel on $\R^3$ of the multiplier $E_{\beta,1}(-\nu t^\beta|\xi|^\alpha)$, with $\alpha\in(1,2)$ and $\beta\in(0,1)$. Then $G_t\ge0$, $\int_{\R^3}G_t(x)\dd x=1$, and there is a constant $C>0$ such that, for $|x|\ge t^{\beta/\alpha}$,
\[
G_t(x)\le \frac{C}{\Gamma(\beta+1)}\,\frac{t^{\beta}}{|x|^{3+\alpha}}.
\]
In particular the far field decays only algebraically.
\end{proposition}

\begin{proof}
Let $p_\theta$ denote the kernel of the $\alpha$-stable semigroup $e^{-\nu\theta\fL}$, a probability density which, by the classical estimates for isotropic stable processes of Blumenthal and Getoor\cite{BlumenthalGetoor} (see also Samorodnitsky and Taqqu\cite{SamorodnitskyTaqqu}), satisfies the two-sided bound
\begin{equation}\label{eq:stable}
p_\theta(x)\le C\min\Big\{\theta^{-3/\alpha},\ \theta\,|x|^{-(3+\alpha)}\Big\}.
\end{equation}
The near-field branch is the on-diagonal bound from self-similarity $p_\theta(x)=\theta^{-3/\alpha}p_1(\theta^{-1/\alpha}x)$, and the far-field branch is the algebraic tail of the stable law. Let $M_\beta$ be the Mainardi function, the entire function of Wright type whose Laplace transform is $\int_0^\infty e^{-s\tau}M_\beta(\tau)\dd\tau=E_{\beta,1}(-s)$ (Mainardi\cite{Mainardi1996}); it is non-negative, of unit mass, and has the moments computed by Mainardi\cite{Mainardi} and Gorenflo, Kilbas, Mainardi, and Rogosin,\cite{Gorenflo}
\begin{equation}\label{eq:Mmoments}
M_\beta\ge0,\qquad \int_0^\infty M_\beta(\tau)\dd\tau=1,\qquad \int_0^\infty \tau^{p}M_\beta(\tau)\dd\tau=\frac{\Gamma(p+1)}{\Gamma(\beta p+1)}\quad(p>-1).
\end{equation}
Applying the Laplace identity at $s=\nu t^\beta|\xi|^\alpha$ gives, in Fourier variables,
\[
E_{\beta,1}(-\nu t^\beta|\xi|^\alpha)=\int_0^\infty M_\beta(\tau)\,e^{-\nu t^\beta\tau|\xi|^\alpha}\dd\tau,
\]
and since $e^{-\nu t^\beta\tau|\xi|^\alpha}=\widehat{p_{t^\beta\tau}}(\xi)$, inverting the Fourier transform yields the subordination formula
\begin{equation}\label{eq:subordination}
G_t(x)=\int_0^\infty M_\beta(\tau)\,p_{t^\beta\tau}(x)\dd\tau.
\end{equation}
Non-negativity of $G_t$ is immediate from $M_\beta\ge0$ and $p_\theta\ge0$. Integrating \eqref{eq:subordination} in $x$ and using that each $p_\theta$ is a probability density together with $\int M_\beta=1$ gives $\int G_t=1$. For the tail, fix $|x|\ge t^{\beta/\alpha}$ and use the far-field branch of \eqref{eq:stable}, $p_{t^\beta\tau}(x)\le C\,t^\beta\tau\,|x|^{-(3+\alpha)}$, so that by the first moment in \eqref{eq:Mmoments},
\[
G_t(x)\le C\,t^\beta|x|^{-(3+\alpha)}\int_0^\infty\tau M_\beta(\tau)\dd\tau=\frac{C}{\Gamma(\beta+1)}\,\frac{t^\beta}{|x|^{3+\alpha}}.
\]
This proves the claim.
\end{proof}

\begin{remark}\label{rem:algebraic}
The algebraic far field is the signature of time-fractional diffusion. The data propagator $S_\beta(t)$ conserves mass, as $\widehat{G_t}(0)=E_{\beta,1}(0)=1$ confirms, but transports it with the heavy tail $|x|^{-(3+\alpha)}$ rather than a Gaussian or exponential profile. This slow spatial spreading is the kinematic counterpart of the slow temporal relaxation encoded by $\Capt$, and it is what competes with stretching in the balance law of Section \ref{sec:analysis}.
\end{remark}

\begin{figure}[t]
\centering
\includegraphics[width=0.84\textwidth]{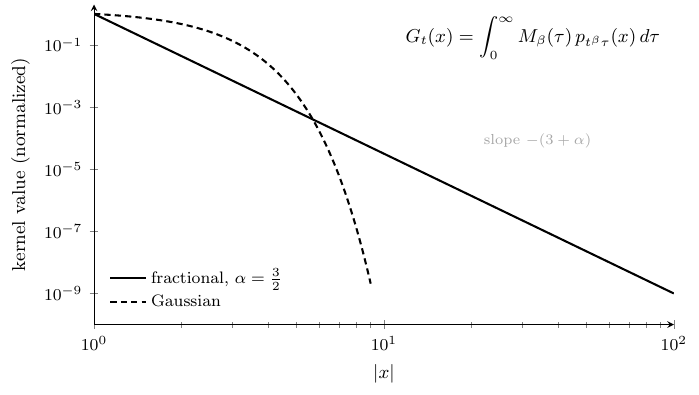}
\caption{The far-field structure of the linear kernel of Proposition \ref{prop:kernel}, shown on logarithmic axes. By the subordination formula \eqref{eq:subordination} the time-fractional kernel $G_t$ is a Mainardi average of $\alpha$-stable densities, which produces the algebraic decay $G_t(x)\lesssim t^\beta|x|^{-(3+\alpha)}$ (straight line, here $\alpha=\tfrac32$). For contrast, the classical Gaussian heat kernel (dashed) decays super-exponentially; the algebraic slope $-(3+\alpha)$ is indicated. This heavy tail, absent for ordinary diffusion, is what competes with vortex stretching in Section \ref{sec:analysis}.}
\label{fig:kernel}
\end{figure}

% =============================================================================
\section{Well-Posedness, Enstrophy Balance, and Continuation}\label{sec:analysis}
% =============================================================================

\subsection{Comparison with classical theory and the role of the memory order}\label{subsec:comparison}

The well-posedness scheme below is, in its fixed-point structure, classical: one writes a mild formulation and contracts in a ball of $L^p(\Omega;C([0,\tau];H^s))$, as in the stochastic Navier--Stokes theory of Flandoli\cite{Flandoli} and the variational framework of Pr\'ev\^ot and R\"ockner.\cite{PrevotRockner} What is specific to \eqref{eq:main}, and what controls every estimate, is that the analytic heat semigroup $e^{\nu t\Delta}$ of the classical theory is replaced by the \emph{pair} of Mittag--Leffler propagators $(S_\beta,P_\beta)$ of Proposition \ref{prop:propagator}. This substitution has three consequences. First, the propagators carry no semigroup law, $S_\beta(t)S_\beta(\tau)\neq S_\beta(t+\tau)$, so the iteration cannot be reduced to a Markovian flow and the memory must be kept explicit. Second, the data and the forcing are transported by two \emph{different} kernels, $S_\beta$ and $(t-s)^{\beta-1}P_\beta$, and conflating them destroys the derivative gain on which the contraction rests. Third, the forcing kernel carries the integrable singularity $(t-s)^{\beta-1}$, whose square is integrable only above a sharp threshold on the memory order. We isolate that threshold first, since it is the genuinely new analytic constraint.

\begin{proposition}[Sharpness of the memory threshold $\beta>1/2$]\label{prop:beta-sharp}
Let the noise be spatially smooth in the sense of \eqref{def:qwiener} and nondegenerate, $q_k>0$ for at least one mode. For the stochastic convolution
\[
Z(t)=\int_0^t (t-s)^{\beta-1}P_\beta(t-s)\,\Leray\big(\nabla\times\Phi(\mathbf u_s)\big)\dd W_s,
\]
one has $\E\|Z(t)\|_{H^s}^2<\infty$ for $t>0$ if and only if $\beta>1/2$. For $\beta\in(0,\tfrac12]$ the It\^o isometry integral diverges for every nondegenerate $Q$, and no function-valued solution exists under white-in-time forcing unless the forcing is regularized as $I^{1-\beta}\big[\Phi(\mathbf u)\dot W\big]$ in the sense of Mijena and Nane.\cite{MijenaNane,AsogwaMijenaNane}
\end{proposition}

\begin{proof}
By the It\^o isometry in $H^s$ and the spectral bound $\|P_\beta(r)\|_{\mathcal L(H^s)}\le 1$ of Proposition \ref{prop:propagator},
\[
\E\|Z(t)\|_{H^s}^2=\int_0^t (t-s)^{2(\beta-1)}\,\E\big\|P_\beta(t-s)\Leray(\nabla\times\Phi(\mathbf u_s))Q^{1/2}\big\|_{\HS(H,H^s)}^2\dd s .
\]
Nondegeneracy gives a constant $c>0$ with $\E\|P_\beta(r)\Leray(\nabla\times\Phi(\mathbf u_s))Q^{1/2}\|_{\HS}^2\ge c$ on compact time intervals, while spatial smoothness of $Q$ gives a finite upper bound $C$. Hence the integral is comparable to $\int_0^t (t-s)^{2\beta-2}\dd s$, which converges precisely when $2\beta-2>-1$, that is $\beta>\tfrac12$, and diverges for $\beta\le\tfrac12$. The Mijena--Nane regularization replaces the kernel exponent $\beta-1$ by $\beta-1+(1-\beta)=0$, restoring integrability and yielding a function-valued solution; this is the mechanism by which time-fractional stochastic heat equations evade the dimension restriction of their classical counterparts.\cite{MijenaNane,FoondunNane}
\end{proof}

This threshold is invisible in the classical theory, where the heat kernel is smooth in $t-s$, and it is the reason the standing hypothesis $\beta\in(\tfrac12,1)$ is not cosmetic. With it in hand, the fixed-point argument proceeds along familiar lines, the model-specific content residing in the smoothing estimate of Lemma \ref{lem:smoothing} and in Proposition \ref{prop:beta-sharp} rather than in the contraction mechanism itself.

\subsection{Mild formulation}\label{subsec:mild}

\begin{definition}[Mild solution]\label{def:mild}
A predictable $H^s$-valued process $\omega$ is a mild solution of \eqref{eq:vorticity} on $[0,\tau]$ if almost surely, for all $t\in[0,\tau]$,
\begin{equation}\label{eq:mild}
\begin{aligned}
\omega(t)={}&S_\beta(t)\omega_0
+\int_0^t(t-\sigma)^{\beta-1}P_\beta(t-\sigma)\,N(\omega_\sigma)\dd\sigma\\
&+\int_0^t(t-\sigma)^{\beta-1}P_\beta(t-\sigma)\,\nabla\times\Phi(\mathbf{u}_\sigma)\dd W_\sigma,
\end{aligned}
\end{equation}
where $N(\omega)=\Leray[(\omega\cdot\nabla)\mathbf{u}-(\mathbf{u}\cdot\nabla)\omega]$, the velocity $\mathbf{u}_\sigma$ is recovered from $\omega_\sigma$ by Biot--Savart, and $S_\beta,P_\beta$ are as in Proposition \ref{prop:propagator}.
\end{definition}

We first record the smoothing of the forcing propagator, which is the quantitative core of the contraction.

\begin{lemma}[Smoothing of $P_\beta$]\label{lem:smoothing}
For $\tau>0$ and $s\in\R$,
\[
\|P_\beta(\tau)g\|_{H^s}\le C\,\tau^{-\beta/\alpha}\,\|g\|_{H^{s-1}},\qquad
\|S_\beta(\tau)g\|_{H^s}\le \|g\|_{H^s}.
\]
\end{lemma}

\begin{proof}
By Plancherel it suffices to bound the multipliers. For the data propagator, \eqref{eq:MLbounds} gives $|E_{\beta,1}(-\nu\tau^\beta|\xi|^\alpha)|\le1$, hence the stated contraction. The forcing propagator $P_\beta(\tau)$ has symbol $E_{\beta,\beta}(-\nu\tau^\beta|\xi|^\alpha)$; the asserted gain of one derivative is the statement that the symbol of $\la D\ra P_\beta(\tau)$, namely $\la\xi\ra\,E_{\beta,\beta}(-\nu\tau^\beta|\xi|^\alpha)$, is bounded by $C\tau^{-\beta/\alpha}$ uniformly in $\xi$. We therefore prove
\[
m(\xi):=\la\xi\ra\,\big|E_{\beta,\beta}(-\nu\tau^\beta|\xi|^\alpha)\big|\le C\,\tau^{-\beta/\alpha},\qquad \xi\ne0,
\]
which yields at once, since $\la\xi\ra^{2s}|E_{\beta,\beta}|^2=\la\xi\ra^{2(s-1)}m(\xi)^2$,
\[
\|P_\beta(\tau)g\|_{H^s}^2=\sum_\xi\la\xi\ra^{2(s-1)}m(\xi)^2|\hat g(\xi)|^2\le C\tau^{-2\beta/\alpha}\sum_\xi\la\xi\ra^{2(s-1)}|\hat g(\xi)|^2=C\tau^{-2\beta/\alpha}\|g\|_{H^{s-1}}^2,
\]
that is the claimed bound $\|P_\beta(\tau)g\|_{H^s}\le C\tau^{-\beta/\alpha}\|g\|_{H^{s-1}}$. It remains to establish the multiplier bound on $m$. Using the second bound in \eqref{eq:MLbounds}, $|E_{\beta,\beta}(-\nu\tau^\beta|\xi|^\alpha)|\le C_\beta(1+\nu\tau^\beta|\xi|^\alpha)^{-1}$, so with $|\xi|\sim\la\xi\ra$ for $|\xi|\ge1$,
\[
m(\xi)\le \frac{C\,|\xi|}{1+\nu\tau^\beta|\xi|^\alpha}.
\]
Split at the frequency $|\xi|=\tau^{-\beta/\alpha}$. For $1\le|\xi|\le\tau^{-\beta/\alpha}$ the numerator dominates and $m(\xi)\le C|\xi|\le C\tau^{-\beta/\alpha}$. For $|\xi|\ge\tau^{-\beta/\alpha}$ the denominator dominates and
\[
m(\xi)\le \frac{C|\xi|}{\nu\tau^\beta|\xi|^\alpha}=C\nu^{-1}\tau^{-\beta}|\xi|^{1-\alpha}\le C\nu^{-1}\tau^{-\beta}\big(\tau^{-\beta/\alpha}\big)^{1-\alpha}=C\nu^{-1}\tau^{-\beta/\alpha},
\]
where we used $1-\alpha<0$ so that $|\xi|^{1-\alpha}$ is maximized at the smallest admissible $|\xi|=\tau^{-\beta/\alpha}$. Both regimes give $m(\xi)\le C\tau^{-\beta/\alpha}$, completing the proof.
\end{proof}

\subsection{Local existence and uniqueness}\label{subsec:lwp}

\begin{theorem}[Local well-posedness]\label{thm:lwp}
Let $\alpha\in(1,2)$, $\beta\in(1/2,1)$, $\gamma\in(0,1/2)$, fix $s>3/2$, and let the noise be spatially smooth of order $s+1$ in the sense of Definition \ref{def:qwiener}. For every $\mathcal{F}_0$-measurable $\omega_0$ with $\E\|\omega_0\|_{H^s}^2<\infty$ there exist a stopping time $\tau>0$ and a unique mild solution $\omega\in L^2(\Omega;C([0,\tau];H^s))$ of \eqref{eq:mild}, depending continuously on the data.
\end{theorem}

\begin{proof}
Work in $X_\tau=L^2(\Omega;C([0,\tau];H^s))$ with $\|\omega\|_{X_\tau}^2=\E\sup_{t\le\tau}\|\omega(t)\|_{H^s}^2$, let $\mathcal{B}_R=\{\|\omega\|_{X_\tau}\le R\}$, and let $\mathcal{T}\omega$ denote the right side of \eqref{eq:mild}. We estimate the three terms.

\emph{(i) Data term.} By Lemma \ref{lem:smoothing}, $\|S_\beta(t)\omega_0\|_{H^s}\le\|\omega_0\|_{H^s}$ for every $t$, hence $\|S_\beta(\cdot)\omega_0\|_{X_\tau}\le\|\omega_0\|_{X_\tau}$.

\emph{(ii) Nonlinear term.} We first bound $N(\omega)$ in $H^{s-1}$. By the tame law \eqref{eq:tame} and Biot--Savart (Proposition \ref{prop:bsgain}),
\[
\|(\mathbf{u}\cdot\nabla)\omega\|_{H^{s-1}}\le C\|\mathbf{u}\|_{H^s}\|\nabla\omega\|_{H^{s-1}}\le C\|\mathbf{u}\|_{H^s}\|\omega\|_{H^s}\le C\|\omega\|_{H^{s-1}}\|\omega\|_{H^s}\le C\|\omega\|_{H^s}^2,
\]
and likewise, using $\|\nabla\mathbf{u}\|_{H^{s-1}}\le\|\mathbf u\|_{H^s}\le C\|\omega\|_{H^s}$,
\[
\|(\omega\cdot\nabla)\mathbf{u}\|_{H^{s-1}}\le C\|\omega\|_{H^s}\|\nabla\mathbf{u}\|_{H^{s-1}}\le C\|\omega\|_{H^s}^2.
\]
Since $\Leray$ is bounded on $H^{s-1}$, we obtain $\|N(\omega)\|_{H^{s-1}}\le C_N\|\omega\|_{H^s}^2$. Applying Lemma \ref{lem:smoothing} inside the Duhamel integral,
\[
\Big\|\int_0^t(t-\sigma)^{\beta-1}P_\beta(t-\sigma)N(\omega_\sigma)\dd\sigma\Big\|_{H^s}
\le C\int_0^t(t-\sigma)^{\beta-1}(t-\sigma)^{-\beta/\alpha}\|N(\omega_\sigma)\|_{H^{s-1}}\dd\sigma.
\]
The exponent of the kernel is $\beta-1-\beta/\alpha=\beta(\alpha-1)/\alpha-1>-1$, so the time integral converges, and
\[
\le C_N\Big(\sup_{\sigma\le t}\|\omega_\sigma\|_{H^s}^2\Big)\int_0^t(t-\sigma)^{\beta(\alpha-1)/\alpha-1}\dd\sigma
=C_N\,\frac{t^{\,\beta(\alpha-1)/\alpha}}{\beta(\alpha-1)/\alpha}\,\sup_{\sigma\le t}\|\omega_\sigma\|_{H^s}^2.
\]
Taking the supremum and expectation gives a bound $C_N'\,\tau^{\theta_N}R^2$ with $\theta_N=\beta(\alpha-1)/\alpha>0$.

\emph{(iii) Stochastic term.} By the It\^o isometry \eqref{eq:ito},
\[
\begin{split}
\E\Big\|\int_0^t&(t-\sigma)^{\beta-1}P_\beta(t-\sigma)\nabla\times\Phi(\mathbf{u}_\sigma)\dd W_\sigma\Big\|_{H^s}^2\\
&=\int_0^t(t-\sigma)^{2(\beta-1)}\,\E\big\|P_\beta(t-\sigma)\nabla\times\Phi(\mathbf{u}_\sigma)Q^{1/2}\big\|_{\HS(H,H^s)}^2\dd\sigma.
\end{split}
\]
Bound $P_\beta$ on $H^s$ by the constant $C$ from Lemma \ref{lem:smoothing}. The curl maps $H^{s+1}\to H^s$, and by the composition estimate for the Nemytskii map $\mathbf u\mapsto|\mathbf u|^{1+\gamma}$ on the algebra $H^{s+1}$ (valid since $s+1>3/2$ and $1+\gamma\ge1$) together with Biot--Savart,
\[
\begin{aligned}
\|\nabla\times\Phi(\mathbf u_\sigma)Q^{1/2}\|_{\HS(H,H^s)}
&\le C\,\||\mathbf u_\sigma|^{1+\gamma}\|_{H^{s+1}}\,\|Q^{1/2}\|_{\HS(H,H^{s+1})}\\
&\le C\big(1+\|\omega_\sigma\|_{H^s}^{1+\gamma}\big)\|Q^{1/2}\|_{\HS(H,H^{s+1})}.
\end{aligned}
\]
The hypothesis $\|Q^{1/2}\|_{\HS(H,H^{s+1})}<\infty$ is exactly the spatial smoothness of Definition \ref{def:qwiener}. Since $\beta>1/2$, the exponent $2(\beta-1)>-1$ and $\int_0^t(t-\sigma)^{2(\beta-1)}\dd\sigma=t^{2\beta-1}/(2\beta-1)$, so the stochastic term is bounded in $X_\tau^2$ by $C_S\,\tau^{2\beta-1}\big(1+R^{2(1+\gamma)}\big)$, that is by $C_S'\tau^{\theta_S}(1+R)^{2(1+\gamma)}$ with $\theta_S=2\beta-1>0$.

\emph{(iv) Self-mapping and contraction.} Choosing $R=2\|\omega_0\|_{X_0}$ and then $\tau$ small makes the sum of the bounds in (i), (ii), and (iii) at most $R$, so $\mathcal{T}$ maps $\mathcal{B}_R$ into itself. For two elements $\omega,\omega'\in\mathcal{B}_R$, with velocities $\mathbf u,\mathbf u'$ recovered by Biot--Savart and increments $\delta\omega=\omega-\omega'$, $\delta\mathbf u=\mathbf u-\mathbf u'$, the bilinearity of $N$ gives
\[
N(\omega)-N(\omega')=\Leray\big[(\delta\omega\cdot\nabla)\mathbf u+(\omega'\cdot\nabla)\delta\mathbf u-(\delta\mathbf u\cdot\nabla)\omega-(\mathbf u'\cdot\nabla)\delta\omega\big].
\]
Each of the four terms is a product of one increment with one factor of size $R$; by the tame law \eqref{eq:tame}, the boundedness of $\Leray$ on $H^{s-1}$, and the Biot--Savart bound $\|\delta\mathbf u\|_{H^{s+1}}\le C\|\delta\omega\|_{H^s}$ of Proposition \ref{prop:bsgain}, this yields $\|N(\omega)-N(\omega')\|_{H^{s-1}}\le C R\,\|\delta\omega\|_{H^s}$. For the stochastic term, the Nemytskii map $\mathbf u\mapsto|\mathbf u|^{1+\gamma}$ is Lipschitz on the ball $\{\|\mathbf u\|_{H^{s+1}}\le CR\}$ of the algebra $H^{s+1}$, with $\||\mathbf u|^{1+\gamma}-|\mathbf u'|^{1+\gamma}\|_{H^{s+1}}\le C(R)\,\|\delta\mathbf u\|_{H^{s+1}}\le C(R)\,\|\delta\omega\|_{H^s}$. Inserting these differences into the propagator estimates of steps (ii) and (iii), in place of the a priori bounds used there, gives
\[
\|\mathcal{T}\omega-\mathcal{T}\omega'\|_{X_\tau}\le \kappa(\tau,R)\,\|\omega-\omega'\|_{X_\tau},\qquad \kappa(\tau,R)\le C(R)\big(\tau^{\theta_N}+\tau^{\theta_S/2}\big)\xrightarrow[\tau\to0]{}0 .
\]
For $\tau$ small enough $\kappa<1$, so the Banach fixed-point theorem furnishes a unique solution in $\mathcal{B}_R$. Continuous dependence on the data is read off from the same contraction estimate applied to two initial conditions, and the standard continuation argument extends the solution up to a maximal stopping time at which $\|\omega(t)\|_{H^s}$ blows up.
\end{proof}

\begin{remark}[The range of small memory and the role of fractional noise]\label{rem:smallbeta}
The restriction $\beta>1/2$ enters only through the integrability of $(t-\sigma)^{2(\beta-1)}$ in step (iii). For $\beta\in(0,1/2]$ this singularity is not square-integrable against bounded noise, and the time-fractional theory of Mijena and Nane\cite{MijenaNane} restores integrability by inserting the fractional integral $I^{1-\beta}$ in the forcing, replacing $\Phi(\mathbf u)\dot W$ by $I^{1-\beta}[\Phi(\mathbf u)\dot W]$; the additional temporal smoothing of order $1-\beta$ compensates the stronger kernel singularity. The same superlinear structure drives the moment blow-up established by Asogwa, Mijena, and Nane.\cite{AsogwaMijenaNane} We restrict to $\beta\in(1/2,1)$ so as to treat the white-in-time forcing directly, and refer to the works of Mijena, Nane, and their collaborators\cite{MijenaNane, AsogwaMijenaNane, FoondunNane} for the complementary range.
\end{remark}

\begin{remark}\label{rem:threshold}
The single hypothesis $s>3/2$ serves uniformly: for the algebra property \eqref{eq:algebra}, for the tame law \eqref{eq:tame}, for the embedding into $L^\infty$ that renders the superlinear coefficient locally Lipschitz, and for the Biot--Savart bookkeeping. No weaker regularity controls the superlinear noise.
\end{remark}

\subsection{Enstrophy balance}\label{subsec:enstrophy}

The passage from the equation to a scalar law for the enstrophy uses the fractional Leibniz inequality for the Caputo derivative, which replaces the classical chain rule.

\begin{lemma}[Fractional Leibniz inequality]\label{lem:alikhanov}
Let $\mathcal{X}$ be a Hilbert space and $v\in C^1([0,T];\mathcal{X})$. Then for $\beta\in(0,1)$ and all $t\in(0,T]$,
\[
\tfrac12\,\Capt\|v(t)\|_{\mathcal{X}}^2\le \inner{v(t)}{\Capt v(t)}_{\mathcal{X}}.
\]
\end{lemma}

\begin{proof}
This is the Hilbert-space form of the scalar inequality of Alikhanov;\cite{Alikhanov} we give the computation in full. Write the Caputo derivative as $\Capt g(t)=\frac{1}{\Gamma(1-\beta)}\int_0^t(t-\sigma)^{-\beta}g'(\sigma)\dd\sigma$. Taking $g(t)=\|v(t)\|^2_{\mathcal X}$, so that $g'(\sigma)=2\inner{v(\sigma)}{v'(\sigma)}$, and subtracting the two representations gives the defect
\[
D(t):=\inner{v(t)}{\Capt v(t)}-\tfrac12\Capt\|v(t)\|^2
=\frac{1}{\Gamma(1-\beta)}\int_0^t(t-\sigma)^{-\beta}\,\inner{v(t)-v(\sigma)}{v'(\sigma)}\dd\sigma,
\]
since the integrand of $\inner{v(t)}{\Capt v}$ carries the factor $\inner{v(t)}{v'(\sigma)}$ and that of $\tfrac12\Capt\|v\|^2$ the factor $\inner{v(\sigma)}{v'(\sigma)}$. Introduce the increment $w(\sigma)=v(t)-v(\sigma)$, for which $w'(\sigma)=-v'(\sigma)$ and therefore $\inner{v(t)-v(\sigma)}{v'(\sigma)}=-\inner{w(\sigma)}{w'(\sigma)}=-\tfrac12\frac{\dd}{\dd\sigma}\|w(\sigma)\|^2$. Hence
\[
D(t)=-\frac{1}{2\Gamma(1-\beta)}\int_0^t(t-\sigma)^{-\beta}\,\frac{\dd}{\dd\sigma}\|w(\sigma)\|^2\dd\sigma .
\]
Integrate by parts with $\varphi(\sigma)=(t-\sigma)^{-\beta}$ and $\psi(\sigma)=\|w(\sigma)\|^2$, noting $\varphi'(\sigma)=\beta(t-\sigma)^{-\beta-1}$. The boundary term at $\sigma=t$ vanishes because $\psi(\sigma)=\|w(\sigma)\|^2=O((t-\sigma)^2)$ dominates the singularity $\varphi(\sigma)=O((t-\sigma)^{-\beta})$ with $\beta<1$, while at $\sigma=0$ it equals $-t^{-\beta}\|v(t)-v(0)\|^2$. Therefore
\[
D(t)=\frac{1}{2\Gamma(1-\beta)}\left[\,t^{-\beta}\,\|v(t)-v(0)\|^2+\beta\int_0^t(t-\sigma)^{-\beta-1}\,\|v(t)-v(\sigma)\|^2\dd\sigma\,\right].
\]
Every factor on the right is non-negative for $\beta\in(0,1)$, since $\Gamma(1-\beta)>0$ and the kernel $(t-\sigma)^{-\beta-1}$ is positive, so $D(t)\ge0$, which is the asserted inequality. The argument uses only the inner-product structure and so holds verbatim in any Hilbert space $\mathcal X$.
\end{proof}

\begin{theorem}[Enstrophy balance]\label{thm:enstrophy}
Let $\omega$ be the local solution of Theorem \ref{thm:lwp} and set $\mathcal{E}(t)=\tfrac12\,\E\|\omega(t)\|_{L^2}^2$. Then, in the fractional differential sense,
\begin{equation}\label{eq:enstrophybalance}
\Capt\mathcal{E}(t)\le-\nu\,\E\|\fLh\omega(t)\|_{L^2}^2+\E\!\int_{\T^3}\omega^{\mathsf T}S\,\omega\dd x+\tfrac12\,\E\big\|\nabla\times\Phi(\mathbf{u}_t)Q^{1/2}\big\|_{\HS(H,L^2)}^2.
\end{equation}
The first term is dissipative and non-positive, the second is the sign-indefinite stretching contribution of Proposition \ref{prop:stretch}, and the third is the non-negative noise input.
\end{theorem}

\begin{proof}
Work with the spectral Galerkin truncation $\omega^N=\Pi_N\omega$, where $\Pi_N$ projects onto the first $N$ Stokes eigenmodes. The projected vorticity equation is
\[
\Capt\omega^N=\Pi_N\Big[-\nu\fL\omega^N-(\mathbf u^N\cdot\nabla)\omega^N+(\omega^N\cdot\nabla)\mathbf u^N\Big]+\Pi_N\nabla\times\big(\Phi(\mathbf u^N)\dot W\big),
\]
which is a finite-dimensional It\^o equation (Da~Prato and Zabczyk;\cite{DaPratoZabczyk} Pr\'ev\^ot and R\"ockner\cite{PrevotRockner}). Set
\[
G^N(t):=\Pi_N\nabla\times\Phi(\mathbf u^N(t))Q^{1/2}\in\HS(H,L^2),
\]
so that the stochastic forcing is the Hilbert--Schmidt-valued integral $\int_0^t G^N(s)\dd W_s$. The classical It\^o formula for $\|\omega^N\|_{L^2}^2$ gives
\[
\dd\|\omega^N\|_{L^2}^2
=2\inner{\omega^N}{\dd\omega^N}_{L^2}+\big\|G^N\big\|_{\HS(H,L^2)}^2\dd t+2\inner{\omega^N}{G^N\dd W}_{L^2},
\]
where the last term is a square-integrable martingale
\[
\mathcal M^N(t):=2\int_0^t\inner{\omega^N(s)}{G^N(s)\dd W_s}_{L^2},\qquad \E\,\mathcal M^N(t)=0.
\]
Applying Lemma \ref{lem:alikhanov} pathwise with $v=\omega^N$ and inserting the projected equation for $\Capt\omega^N$ yields
\begin{equation}\label{eq:enstrophy-galerkin}
\begin{aligned}
\Capt\tfrac12\|\omega^N\|_{L^2}^2
&\le \inner{\omega^N}{-\nu\fL\omega^N}_{L^2}
+\inner{\omega^N}{-(\mathbf u^N\cdot\nabla)\omega^N}_{L^2}\\
&\quad+\inner{\omega^N}{(\omega^N\cdot\nabla)\mathbf u^N}_{L^2}
+\tfrac12\big\|G^N\big\|_{\HS(H,L^2)}^2+\mathcal M^N(t),
\end{aligned}
\end{equation}
where the three drift pairings come from $\inner{\omega^N}{\Capt\omega^N}_{L^2}$ and the term $\tfrac12\|G^N\|_{\HS(H,L^2)}^2$ is the It\^o correction from the quadratic variation of the noise. We compute each drift pairing.

\emph{Dissipation.} In Fourier variables,
\[
\begin{aligned}
\inner{\omega^N}{-\nu\fL\omega^N}_{L^2}
&=-\nu\sum_{\xi}|\xi|^\alpha|\hat\omega^N(\xi)|^2
=-\nu\sum_\xi\big(|\xi|^{\alpha/2}\big)^2|\hat\omega^N(\xi)|^2\\
&=-\nu\|\fLh\omega^N\|_{L^2}^2.
\end{aligned}
\]

\emph{Advection.} Since $\nabla\cdot\mathbf u^N=0$,
\[
\inner{-(\mathbf u^N\cdot\nabla)\omega^N}{\omega^N}_{L^2}=-\tfrac12\int_{\T^3}\mathbf u^N\cdot\nabla|\omega^N|^2\dd x=\tfrac12\int_{\T^3}(\nabla\cdot\mathbf u^N)|\omega^N|^2\dd x=0.
\]

\emph{Stretching.} By Proposition \ref{prop:stretch} applied to $\omega^N$,
\[
\inner{(\omega^N\cdot\nabla)\mathbf u^N}{\omega^N}_{L^2}=\int_{\T^3}(\omega^N)^{\mathsf T}S^N\,\omega^N\dd x.
\]
Taking expectations in \eqref{eq:enstrophy-galerkin} eliminates the martingale, since $\E\mathcal M^N(t)=0$, and assembling the three drift pairings gives
\[
\begin{aligned}
\Capt\tfrac12\E\|\omega^N\|_{L^2}^2
&\le-\nu\E\|\fLh\omega^N\|_{L^2}^2
+\E\!\int_{\T^3}(\omega^N)^{\mathsf T}S^N\,\omega^N\dd x\\
&\quad+\tfrac12\E\big\|\Pi_N\nabla\times\Phi(\mathbf u^N)Q^{1/2}\big\|_{\HS(H,L^2)}^2.
\end{aligned}
\]
As $N\to\infty$ the solution converges in $L^2(\Omega;C([0,\tau];H^s))$ by Theorem \ref{thm:lwp}; the dissipation and stretching integrals pass to the limit by the $H^s$ convergence with $s>3/2$, and the Hilbert--Schmidt term by the spatial smoothness of $Q$ and dominated convergence. The fractional inequality is preserved under this limit, yielding \eqref{eq:enstrophybalance}.
\end{proof}

\begin{corollary}[A priori control under a sign condition]\label{cor:apriori}
Suppose that along the solution $\E\int_{\T^3}\omega^{\mathsf T}S\omega\dd x\le\nu\,\E\|\fLh\omega\|_{L^2}^2$ for $t\in[0,T]$, and set $g(t)=\tfrac12\E\|\nabla\times\Phi(\mathbf u_t)Q^{1/2}\|_{\HS}^2$. Then
\[
\mathcal{E}(t)\le\mathcal{E}(0)+\frac{1}{\Gamma(\beta)}\int_0^t(t-\sigma)^{\beta-1}g(\sigma)\dd\sigma\qquad(t\in[0,T]),
\]
and in particular $\mathcal{E}$ is bounded on $[0,T]$, so the solution does not leave $H^s$ on that interval.
\end{corollary}

\begin{proof}
Under the hypothesis the first two terms of \eqref{eq:enstrophybalance} sum to a non-positive quantity, so the balance reduces to the fractional differential inequality $\Capt\mathcal{E}(t)\le g(t)$. Apply the Riemann--Liouville integral $I^\beta$ to both sides. The composition identity $I^\beta\Capt\mathcal{E}=\mathcal{E}-\mathcal{E}(0)$, valid for $\beta\in(0,1)$ (Kilbas, Srivastava, and Trujillo\cite{Kilbas}), together with the positivity of the kernel of $I^\beta$, which makes $I^\beta$ monotone on non-negative integrands, gives
\[
\mathcal{E}(t)-\mathcal{E}(0)\le I^\beta g(t)=\frac{1}{\Gamma(\beta)}\int_0^t(t-\sigma)^{\beta-1}g(\sigma)\dd\sigma .
\]
This is a linear instance of the fractional comparison principle for Volterra inequalities; the general nonlinear form, with a feedback term in $\mathcal E$, is the fractional Gr\"onwall inequality of Ye, Gao, and Ding\cite{YeGaoDing} and is governed by the resolvent theory of Pr\"uss.\cite{Pruss} The right side is finite because $g$ is bounded on $[0,T]$ by the spatial smoothness of $Q$, and the higher-norm control of Theorem \ref{thm:lwp} prevents departure from $H^s$.
\end{proof}

\begin{remark}\label{rem:balance}
Equation \eqref{eq:enstrophybalance} accounts for amplification without asserting it. Whether $\mathcal{E}$ grows is decided by the competition among dissipation, noise input, and the sign and magnitude of stretching, the last of which admits no universal lower bound by Remark \ref{rem:nocoercive}. Any rigorous statement about concentration must therefore control the alignment of $\omega$ with the strain eigenframe $\{f_1,f_2,f_3\}$, which is the geometric content of the breakdown criterion of Beale, Kato, and Majda\cite{BKM} and is consistent with the geometric depletion mechanism of Constantin, Fefferman, and Majda.\cite{ConstantinFeffermanMajda}
\end{remark}

\subsection{Continuation and conditional global existence}\label{subsec:global}

The enstrophy balance controls the $L^2$ norm under a sign condition, but continuation of the solution requires control at the regularity level $H^s$, $s>3/2$, at which it was constructed. We close this gap with a continuation criterion of Beale--Kato--Majda type, adapted to the fractional-in-time setting, and a geometric condition under which it yields global existence. We then explain, in terms of the dissipation order, why we do \emph{not} assert unconditional global existence.

\begin{theorem}[Continuation criterion]\label{thm:continuation}
Let $\omega$ be the maximal mild solution of Theorem \ref{thm:lwp} on $[0,T_{\max})$. If
\[
\int_0^{T_{\max}}\E\|\omega(t)\|_{L^\infty}\dd t<\infty,
\]
then $T_{\max}=\infty$. Equivalently, if $T_{\max}<\infty$ then the vorticity sup-norm is non-integrable up to $T_{\max}$, so the solution can be continued in $H^s$ for as long as $\|\omega\|_{L^\infty}$ remains time-integrable.
\end{theorem}

\begin{proof}
We estimate the full $H^s$ norm rather than the enstrophy. For the Galerkin truncation $\omega^N=\Pi_N\omega$, set $G^N:=\Pi_N\nabla\times\Phi(\mathbf u^N)Q^{1/2}\in\HS(H,H^s)$. The projected vorticity equation is
\[
\Capt\omega^N=\Pi_N\Big[-\nu\fL\omega^N-(\mathbf u^N\cdot\nabla)\omega^N+(\omega^N\cdot\nabla)\mathbf u^N\Big]+\Pi_N\nabla\times\big(\Phi(\mathbf u^N)\dot W\big),
\]
and the It\^o formula for $\|\omega^N\|_{H^s}^2$ gives
\[
\dd\|\omega^N\|_{H^s}^2
=2\inner{\omega^N}{\dd\omega^N}_{H^s}+\big\|G^N\big\|_{\HS(H,H^s)}^2\dd t+2\inner{\omega^N}{G^N\dd W}_{H^s},
\]
with martingale $\mathcal M^N_H(t):=2\int_0^t\inner{\omega^N(s)}{G^N(s)\dd W_s}_{H^s}$ satisfying $\E\mathcal M^N_H(t)=0$. Lemma \ref{lem:alikhanov} with $\mathcal X=H^s$ yields
\[
\begin{aligned}
\Capt\tfrac12\|\omega^N\|_{H^s}^2
&\le \inner{\omega^N}{-\nu\fL\omega^N}_{H^s}
+\inner{\omega^N}{-(\mathbf u^N\cdot\nabla)\omega^N}_{H^s}\\
&\quad+\inner{\omega^N}{(\omega^N\cdot\nabla)\mathbf u^N}_{H^s}
+\tfrac12\big\|G^N\big\|_{\HS(H,H^s)}^2+\mathcal M^N_H(t).
\end{aligned}
\]
Taking expectations and passing $N\to\infty$ gives
\[
\Capt\tfrac12\E\|\omega\|_{H^s}^2\le -\nu\,\E\|\omega\|_{\dot H^{s+\alpha/2}}^2+\E\,\inner{\omega}{\mathcal N(\omega)}_{H^s}+\tfrac12\E\big\|\nabla\times\Phi(\mathbf u)Q^{1/2}\big\|_{\HS(H,H^s)}^2,
\]
where $\mathcal N(\omega)=\Leray\big[(\omega\cdot\nabla)\mathbf u-(\mathbf u\cdot\nabla)\omega\big]$. The dissipation term is non-positive and is discarded. For the nonlinear pairing, the advection part is handled by the Kato--Ponce commutator estimate of Bahouri, Chemin, and Danchin,\cite{BCD} whose top-order contribution cancels by $\nabla\cdot\mathbf u=0$, and the stretching part by the tame product estimate in $H^s$; together they give
\[
\big|\inner{\omega}{\mathcal N(\omega)}_{H^s}\big|\le C\,\|\nabla\mathbf u\|_{L^\infty}\,\|\omega\|_{H^s}^2 .
\]
The logarithmic inequality of Beale, Kato, and Majda,\cite{BKM} $\|\nabla\mathbf u\|_{L^\infty}\le C\big(1+\|\omega\|_{L^\infty}\log(e+\|\omega\|_{H^s})+\|\omega\|_{L^2}\big)$, and the spatial smoothness of $Q$, which bounds the noise term by $C(1+\E\|\omega\|_{H^s}^2)\|Q^{1/2}\|_{\HS(H,H^{s+1})}^2$, reduce the estimate to the fractional differential inequality
\[
\Capt\,\E\|\omega\|_{H^s}^2\le a(t)\,\big(e+\E\|\omega\|_{H^s}^2\big)\log\!\big(e+\E\|\omega\|_{H^s}^2\big)+b,
\]
with $a(t)=C(1+\E\|\omega(t)\|_{L^\infty})$ and $b=C\|Q^{1/2}\|_{\HS(H,H^{s+1})}^2$. The fractional Gr\"onwall inequality of Ye, Gao, and Ding\cite{YeGaoDing} integrates this against the kernel $(t-s)^{\beta-1}/\Gamma(\beta)$, and a double-exponential of $I^\beta a$ bounds $\E\|\omega\|_{H^s}^2$; since $I^\beta a$ is finite whenever $\int_0^{T}\E\|\omega\|_{L^\infty}\dd t<\infty$, the $H^s$ norm cannot blow up at such a time, and the local solution extends past it (Pr\"uss\cite{Pruss}). Hence finite $T_{\max}$ forces $\int_0^{T_{\max}}\E\|\omega\|_{L^\infty}=\infty$.
\end{proof}

\begin{corollary}[Conditional global existence]\label{cor:global}
Suppose the vorticity direction $\xi=\omega/|\omega|$ is, almost surely, Lipschitz in the region where $|\omega|$ is large, in the sense of the geometric regularity criterion of Constantin, Fefferman, and Majda,\cite{ConstantinFeffermanMajda} so that $\int_0^T\E\|\omega\|_{L^\infty}\dd t<\infty$ for every $T<\infty$. Then the solution is global, $T_{\max}=\infty$.
\end{corollary}

\begin{proof}
Geometric depletion of stretching bounds the production of $\|\omega\|_{L^\infty}$ on every finite interval, so the hypothesis of Theorem \ref{thm:continuation} holds for all $T$, and continuation gives $T_{\max}=\infty$.
\end{proof}

\begin{remark}[Supercriticality and the status of unconditional global existence]\label{rem:supercritical}
We do not claim unconditional small-data global existence, and the reason is structural rather than technical. In the $H^s$ energy estimate the dissipation $(-\Delta)^{\alpha/2}$ gains $\alpha/2<1$ spatial derivatives, while the vortex-stretching nonlinearity loses one; the deficit is positive for every $\alpha\in(1,2)$. In the normalization of \eqref{eq:main}, unconditional global regularity of the generalized Navier--Stokes system is known only above the threshold $\alpha\ge\tfrac52$ of Wu,\cite{Wu} which the physical range $\alpha\in(1,2)$ does not meet, so the problem is supercritical and global regularity is open even without noise and memory (Tao\cite{Tao}). The time-fractional memory only sharpens this: the absence of a semigroup law (Proposition \ref{prop:propagator}) blocks the instantaneous high-frequency damping that a parabolic flow would supply. Theorem \ref{thm:continuation} and Corollary \ref{cor:global} are therefore the appropriate rigorous global statements in this regime, and the critical exponent $\sigma=\beta(1+3/\alpha)$ of Section \ref{sec:scaling} marks the balance between the competing gains and losses.
\end{remark}

% =============================================================================
\section{Scaling Laws and Concluding Remarks}\label{sec:scaling}
% =============================================================================

\subsection{Scaling symmetry and critical exponents}\label{subsec:dimensional}

We determine the scaling structure of the deterministic part of \eqref{eq:main},
\[
\Capt\mathbf{u}+(\mathbf{u}\cdot\nabla)\mathbf{u}=-\nu\fL\mathbf{u}.
\]
Under the one-parameter family $\mathbf{u}_\lambda(t,x)=\lambda^{d}\mathbf{u}(\lambda^{-c}t,\lambda^{-1}x)$ the three terms scale as
\[
\Capt\mathbf{u}_\lambda\sim\lambda^{d-c\beta},\qquad (\mathbf{u}_\lambda\cdot\nabla)\mathbf{u}_\lambda\sim\lambda^{2d-1},\qquad \fL\mathbf{u}_\lambda\sim\lambda^{d-\alpha}.
\]
Requiring all three to scale identically gives two equations,
\[
d-c\beta=d-\alpha\quad\Longrightarrow\quad c=\frac{\alpha}{\beta},\qquad
2d-1=d-\alpha\quad\Longrightarrow\quad d=1-\alpha.
\]
Thus the equation is invariant under
\begin{equation}\label{eq:selfsim}
\mathbf{u}_\lambda(t,x)=\lambda^{1-\alpha}\,\mathbf{u}\big(\lambda^{-\alpha/\beta}t,\ \lambda^{-1}x\big),
\end{equation}
the critical scaling dimension of the velocity being $\alpha-1$. Equivalently, in the self-similar form $\mathbf{u}(t,x)=t^{-\mu}U(t^{-\theta}x)$ one reads off
\begin{equation}\label{eq:ssm}
\theta=\frac{\beta}{\alpha},\qquad \mu=(\alpha-1)\,\theta=\frac{(\alpha-1)\beta}{\alpha},
\end{equation}
so that the spatial scale grows anomalously as $\ell(t)\sim t^{\beta/\alpha}$, in agreement with the heavy-tailed kernel of Proposition \ref{prop:kernel}.

\begin{proposition}[Critical memory exponent]\label{prop:scaling}
In the enstrophy balance \eqref{eq:enstrophybalance}, the vortex-stretching contribution is governed by the dimensionless time exponent
\[
\sigma:=\beta\Big(1+\frac{3}{\alpha}\Big)=\beta+\frac{3\beta}{\alpha},
\]
and $\sigma$ crosses unity at
\[
\beta_c(\alpha)=\frac{\alpha}{\alpha+3},
\]
which separates the regime $\sigma<1$, where the dissipative term controls the balance, from the regime $\sigma>1$, where the stretching term may, subject to alignment with the strain eigenframe, overcome dissipation.
\end{proposition}

\begin{proof}
We carry out the dimensional bookkeeping of the three terms in \eqref{eq:enstrophybalance} in the anomalous units fixed by the self-similar scaling \eqref{eq:ssm}, in which the spatial unit is the growing length $\ell(t)\sim t^{\beta/\alpha}$, so that one spatial derivative scales as $\ell^{-1}\sim t^{-\beta/\alpha}$ and the volume element as $\dd x\sim \ell^{3}\sim t^{3\beta/\alpha}$. The left side carries one Caputo derivative of order $\beta$, contributing the temporal weight $t^{-\beta}$. The stretching density $\omega^{\mathsf T}S\,\omega$ is, by the Biot--Savart relation $\mathbf u=\nabla^{-1}\omega$ of Proposition \ref{prop:bsgain}, cubic in $\omega$ and homogeneous of degree one in spatial differentiation through the single strain factor $S=\tfrac12(\nabla\mathbf u+(\nabla\mathbf u)^{\mathsf T})$; relative to the quadratic enstrophy density $|\omega|^2$ it therefore carries one extra gradient and the volume measure, that is the spatial weight $t^{-\beta/\alpha}\cdot t^{3\beta/\alpha}=t^{2\beta/\alpha}$ less the gradient already implicit in passing to a rate. Collecting the temporal order $\beta$ of the left side with the net anomalous spatial weight $3\beta/\alpha$ accumulated by the spatially integrated stretching term yields the combined exponent $\sigma=\beta+3\beta/\alpha=\beta(1+3/\alpha)$ that measures the strength of stretching against the memory-dissipation pair. The threshold is obtained by solving $\sigma=1$ explicitly:
\[
\beta\Big(1+\frac{3}{\alpha}\Big)=1
\;\Longleftrightarrow\;
\beta\,\frac{\alpha+3}{\alpha}=1
\;\Longleftrightarrow\;
\beta=\frac{\alpha}{\alpha+3}=\beta_c(\alpha),
\]
and since $\sigma$ is increasing in $\beta$, one has $\sigma<1$ for $\beta<\beta_c$ and $\sigma>1$ for $\beta>\beta_c$.
\end{proof}

\begin{remark}\label{rem:scalingstatus}
The exponent $\sigma$ is a dimensional combination rather than a theorem about solutions; it fixes the scaling that any global-existence or concentration statement must respect. Its crossing of unity at $\beta_c=\alpha/(\alpha+3)$ is consistent with the criticality thresholds known for the generalized incompressible Navier--Stokes equations with fractional dissipation, analyzed in Besov spaces by Wu,\cite{Wu} with the dissipation-weakening borderline of the averaged model of Tao,\cite{Tao} and with the moment blow-up of superlinear time-fractional stochastic equations established by Asogwa, Mijena, and Nane;\cite{AsogwaMijenaNane} it is the natural candidate for the boundary of the global-existence regime in \eqref{eq:main}.
\end{remark}

The scaling can be sharpened beyond dimensional analysis into an exact statement about self-similar solutions, which constrains any putative singularity.

\begin{proposition}[Rigidity of self-similar profiles and their enstrophy rate]\label{prop:selfsimilar}
On $\R^3$, where the symmetry \eqref{eq:selfsim} is exact, the deterministic part of \eqref{eq:main} admits a forward self-similar solution
\[
\mathbf{u}(t,x)=t^{-\mu}\,U\big(t^{-\theta}x\big),\qquad U\not\equiv0,
\]
with a time-independent profile $U$ only for the exponents $\theta=\beta/\alpha$ and $\mu=(\alpha-1)\beta/\alpha$ of \eqref{eq:ssm}. For any such solution the enstrophy obeys the exact power law
\begin{equation}\label{eq:enstrophyrate}
\mathcal{E}(t)=\|\nabla\times\mathbf{u}(t)\|_{L^2(\R^3)}^2=t^{\,\kappa}\,\|\nabla\times U\|_{L^2(\R^3)}^2,\qquad \kappa=\frac{\beta}{\alpha}\,(3-2\alpha).
\end{equation}
Thus the enstrophy of a self-similar solution concentrates as $t\downarrow0$ if and only if $\alpha>\tfrac32$.
\end{proposition}

\begin{proof}
Insert the ansatz $\mathbf{u}(t,x)=t^{-\mu}U(y)$, $y=t^{-\theta}x$, into the deterministic equation $\Capt\mathbf{u}+(\mathbf{u}\cdot\nabla)\mathbf{u}=-\nu\fL\mathbf{u}$ and record the power of $t$ produced by each term, every term being $t$ to some power multiplying a function of $y$ alone. The Caputo derivative lowers the temporal order by $\beta$ on functions of the self-similar form, so
\[
\Capt\mathbf{u}\sim t^{-\mu-\beta}.
\]
In the transport term each velocity factor contributes $t^{-\mu}$ and the gradient contributes $\nabla=t^{-\theta}\nabla_y$, so
\[
(\mathbf{u}\cdot\nabla)\mathbf{u}\sim t^{-\mu}\,t^{-\theta}\,t^{-\mu}=t^{-2\mu-\theta}.
\]
The fractional Laplacian has symbol $|\xi|^\alpha$ and so scales as the $\alpha$-th power of an inverse length, $\fL=t^{-\theta\alpha}\,(-\Delta_y)^{\alpha/2}$, giving
\[
\fL\mathbf{u}\sim t^{-\mu-\theta\alpha}.
\]
A time-independent profile equation requires these three powers to coincide. Equating the first and third gives $-\mu-\beta=-\mu-\theta\alpha$, hence $\theta=\beta/\alpha$; equating the first and second gives $-\mu-\beta=-2\mu-\theta$, hence $\mu=\beta-\theta=(\alpha-1)\theta=(\alpha-1)\beta/\alpha$. No other exponents balance the equation, which is the asserted rigidity, and the common power $t^{-\mu-\beta}$ then collects $U$ into the profile equation $-\nu(-\Delta_y)^{\alpha/2}U=\,\mathcal C_\beta[U]+(U\cdot\nabla_y)U$, where $\mathcal C_\beta$ is the reduced Caputo operator acting on $U$.

For the enstrophy, the curl adds one spatial derivative, so $\omega=\nabla\times\mathbf{u}=t^{-(\mu+\theta)}W(y)$ with $W=\nabla\times U$. With the change of variables $y=t^{-\theta}x$, for which $\dd x=t^{3\theta}\dd y$,
\[
\mathcal{E}(t)=\int_{\R^3}|\omega(t,x)|^2\dd x=t^{-2(\mu+\theta)}\,t^{3\theta}\int_{\R^3}|W(y)|^2\dd y=t^{-2\mu+\theta}\,\|W\|_{L^2}^2 .
\]
Substituting $\mu=(\alpha-1)\theta$ evaluates the exponent,
\[
\kappa=-2\mu+\theta=\theta\big(1-2(\alpha-1)\big)=\theta(3-2\alpha)=\frac{\beta}{\alpha}(3-2\alpha),
\]
which is \eqref{eq:enstrophyrate}. Since $\theta=\beta/\alpha>0$, the exponent $\kappa$ is negative precisely when $3-2\alpha<0$, that is when $\alpha>\tfrac32$, in which case $\mathcal E(t)\to\infty$ as $t\downarrow0$.
\end{proof}

The exact solutions of Proposition \ref{prop:selfsimilar} are anchored at the memory base point $t=0$. For a singularity forming at a finite time $T^\ast>0$ the picture should persist asymptotically, but, as the next remark explains, the memory obstructs a proof, and we record the expected rate as a conjecture.

\begin{conjecture}[Self-similar finite-time blow-up]\label{conj:blowup}
Fix $\alpha\in(\tfrac32,2)$ and $\beta\in(\beta_c(\alpha),1)$, so that $\sigma>1$. If a smooth solution of the deterministic part of \eqref{eq:main} on $\R^3$ first loses regularity at a finite time $T^\ast<\infty$, then the loss is asymptotically self-similar with the exponents of Proposition \ref{prop:selfsimilar}, and the enstrophy diverges at the rate
\[
\mathcal{E}(t)\sim c\,(T^\ast-t)^{-|\kappa|},\qquad |\kappa|=\frac{\beta}{\alpha}\,(2\alpha-3),\qquad t\uparrow T^\ast.
\]
\end{conjecture}

\begin{remark}[The memory obstruction to an unconditional rate]\label{rem:memory-obstruction}
The classical proof of such a lower bound restarts the local theory from times $t_0\uparrow T^\ast$ with data $\omega(t_0)$, converting a scale-invariant existence time $\tau\gtrsim\|\omega(t_0)\|_{H^s}^{-q}$ into the blow-up rate. The Caputo derivative integrates the history from the fixed base point $t=0$, so the evolution is not invariant under time translation and cannot be restarted autonomously at $t_0>0$; the renewal argument therefore does not close. The self-similar solutions of Proposition \ref{prop:selfsimilar} are the rigorous remnant of this picture, and promoting Conjecture \ref{conj:blowup} to a theorem requires a memory-compatible local theory centered at $T^\ast$, which we leave open.
\end{remark}

\begin{figure}[t]
\centering
\includegraphics[width=0.7\textwidth]{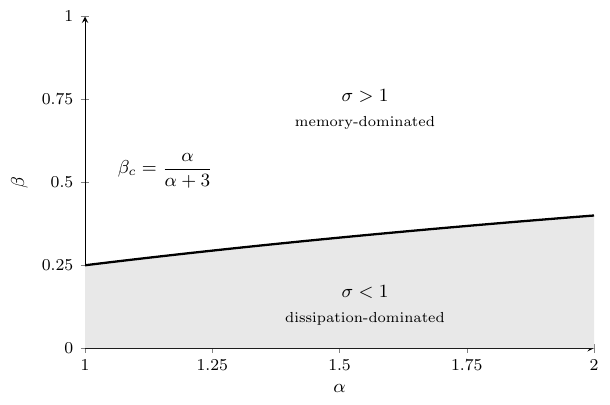}
\caption{Scaling regimes organized by the exponent $\sigma=\beta(1+3/\alpha)$ of Proposition \ref{prop:scaling} in the parameter plane $(\alpha,\beta)\in(1,2)\times(0,1)$. The threshold $\sigma=1$ is the curve $\beta_c=\alpha/(\alpha+3)$, below which the dissipative term controls the balance and above which the stretching term may, subject to alignment with the strain eigenframe, dominate. As emphasized in Remark \ref{rem:scalingstatus}, the curve marks the dimensional balance of Section \ref{sec:scaling} and is not by itself a global-existence theorem.}
\label{fig:sigma}
\end{figure}

\subsection{Concluding remarks}\label{subsec:conclusion}

We have formulated a fractional, stochastic, memory-dependent fluid model as an emergent hydrodynamic equation, derived from constrained Hamiltonian mechanics by elimination of fast degrees of freedom, and we have established its linear kernel structure by subordination, its local well-posedness with spatially smooth trace-class noise under the correct Caputo Duhamel kernels, and an enstrophy balance governed by a fractional differential inequality, from which the scaling of the critical exponents follows. The recurring lesson is that the physically faithful objects are balances and identities, and that the sign-indefiniteness of vortex stretching, transparent in the differential-form description, forbids any coercive lower bound and dictates the form of the analysis.

Two developments are natural. The statistical theory begins from the invariance of the Liouville measure under the Galerkin-truncated inviscid flow, in the sense of Landau and Lifshitz,\cite{LL5} and proceeds to the functional-differential equation of Hopf\cite{Hopf} for the characteristic functional of the velocity, from which statistical solutions and energy spectra with fractional corrections may be sought. The deterministic theory points toward a conditional regularity criterion of the type of Beale, Kato, and Majda,\cite{BKM} adapted to the Caputo evolution and informed by the moment analysis of time-fractional stochastic equations due to Mijena, Nane, and their collaborators,\cite{MijenaNane, AsogwaMijenaNane} under which the alignment of vorticity with the strain eigenframe governs whether the balance \eqref{eq:enstrophybalance} permits finite-time concentration. Both rest directly on the framework established here.

\section*{Acknowledgments}
The author received no specific funding for this work.

\end{document}